\documentclass[preprint,12pt,3p]{elsarticle-no-natbib}

\usepackage{amssymb}
\usepackage{amsthm}

\usepackage[backend=bibtex,%
style=numeric,%
bibencoding=ascii
eprint=false,isbn=false,doi=false, url=false]{biblatex}

\renewbibmacro{in:}{%
   \ifentrytype{article}{}{\printtext{\bibstring{in}\intitlepunct}}}
  
\usepackage{amsmath}
\usepackage{amsfonts}
\usepackage[utf8]{inputenc}
\usepackage[T1]{fontenc}
\usepackage{graphicx}
\usepackage{color}
\usepackage{enumerate}
\usepackage{pdfpages}
\usepackage{harpoon}
\usepackage{subcaption}
\newtheorem{thm}{Theorem}[section]

\newtheorem{lemma}[thm]{Lemma}

\theoremstyle{definition}

\newtheorem{observation}[thm]{Observation}

\theoremstyle{remark}

\numberwithin{equation}{section}
	{\vskip\topsep\noindent\textsc{#1.\/}}%
	{\par\vskip\topsep}%
\definecolor{dkgreen}{rgb}{0,0.7,0.1}
\long\def\ignore#1{}

\def\realline{\hbox to \hsize}
\journal{Symmetry}
\makeatletter
\def\ps@pprintTitleNoSubTo{%
     \let\@oddhead\@empty
     \let\@evenhead\@empty
     \def\@oddfoot{\footnotesize\itshape\hfill\today}%
     \let\@evenfoot\@oddfoot}
\let\ps@pprintTitle\ps@pprintTitleNoSubTo
\makeatother

\def\mZ{\mathbb{Z}}
\def\mR{\mathbb{R}}
\def\cB{\mathcal{B}}
\def\cD{\mathcal{D}}
\def\cM{\mathcal{M}}
\def\cP{\mathcal{P}}
\def\sym{{\fam0 Sym}}

\def\iv{^{-1}}
\let\al\alpha
\let\ma\mu 
\let\sr=q 
\let\Cs\Sigma 
 \def\Csd{\Cs^{\fam0 d}} 
 \def\Csb{\Cs^{\fam0 b}} 
 \def\Csa{\Cs^{\fam0 a}} 
\let\Ga\Gamma 
\let\Wg\Omega 
 \def\Wgd{\Wg^{\fam0 d}} 
 \def\Wgb{\Wg^{\fam0 b}} 
 \def\Wga{\Wg^{\fam0 a}} 
\def\Ad(#1){\delta^{#1}}
\def\Adx(#1){\Ad({#1})(n)}
\def\Ab(#1){\beta^{#1}}
\def\Abx(#1){\Ab({#1})(n,k)}
\def\Aa(#1){\alpha^{#1}}
\def\Aax(#1){\Aa({#1})(n,k)}
\def\wcup{\cup} 
\def\Sz_#1{\sym(\mZ_{#1})}
\let\ga\gamma
\let\col\psi 
\let\si\sigma

\let\eul\varphi
\def\dv{\,|\,}
\let\la\langle
\let\ra\rangle
\def\fix{{\fam0 Fix}}
\long\def\ignore#1{}

\def\pz_#1{ {\mZ_{#1} \times \mZ_{#1}} }
\def\pr_#1{ {\mR_{#1} \times \mR_{#1}} }
\def\symz_#1{ {\sym(\mZ_{#1})} }
\def\symr_#1{ {\sym(\mR_{#1})} }
\def\Ztn{\mZ_{2n}}
\let\lab\lambda
\def\st{\;|\;\allowbreak}
\def\set#1{\{#1\}}

\def\dg#1{\overrightharp{\ensuremath{#1}}}
\let\ep\varepsilon
\def\cN{\mathcal{N}}
\def\cA{\mathcal{A}}

\def\cS{\mathcal{S}}
\let\nsubeq\trianglelefteq
\let\De\Delta

\def\rc#1{\text{\large$\textstyle\frac{1}{#1}$}}
\def\ig#1#2{}

\def\values#1 #2.{\smallskip\noindent{\textbf{Values for $0 \le n \le
12$:\/} #2.}}
\def\valuesk#1 #2.{\smallskip\noindent{\textbf{Values for $k=1$ and $0 \le
n \le 12$:\/} #2.}}

\def\svalues#1 #2.{\smallskip\noindent{\textbf{S values for $0 \le n \le
12$:\/} #2.}}
\def\apvalues#1 #2.{\par\noindent{\textbf{AP values for $0 \le n \le
12$:\/} #2.}}
\def\svaluesk#1 #2.{\smallskip\noindent{\textbf{S values for $k=1$ and $0
\le n \le 12$:\/} #2.}}
\def\apvaluesk#1 #2.{\par\noindent{\textbf{AP values for $k=1$ and $0
\le n \le 12$:\/} #2.}}

\def\cS{\mathcal{S}}
\let\nsubeq\trianglelefteq
\let\De\Delta

\newenvironment{computing}[1]%
        {\removelastskip\vskip\topsep\noindent\textbf{Computing
#1.\/}}%
        {\par\vskip\topsep}%
\newenvironment{counting}[1]%
        {\removelastskip\vskip\topsep\noindent\textbf{Counting #1.\/}}%
        {\par\vskip\topsep}%
        
\newenvironment{countitem}%
       {\removelastskip\vskip\topsep\noindent\ignorespaces}%
       {\par\vskip\topsep}%

\def\nn{n \times n}
\def\tx(#1){\text{(#1)}}

\newenvironment{tightequation}%
    {\begin{equation*}\abovedisplayskip=2pt\belowdisplayskip=2pt%
       \abovedisplayshortskip=2pt\belowdisplayshortskip=2pt}%
    {\end{equation*}}

\def\hp#1#2{\par\noindent\hangindent=#2\parindent\hskip#1\parindent}

\def\listitem{\hp12}

\begin{document}

\begin{frontmatter}

\title{A Catalog of Enumeration Formulas for Bouquet and Dipole Embeddings Under Symmetries}


\author[label1, label3]{M. N. Ellingham\corref{cor1}}

\address[label1]{Department of Mathematics, 1326 Stevenson Center,
Vanderbilt University,
Nashville, Tennessee 37240}
\fntext[label3]{Supported by Simons Foundation award 429625}

\ead{mark.ellingham@vanderbilt.edu} \ead[url]{https://math.vanderbilt.edu/ellingmn/}

\author[label2,label4]{Joanna A. Ellis-Monaghan}

\address[label2]{Korteweg-de Vries Institute for Mathematics, University of Amsterdam, Science Park 105-107, 1098 XH Amsterdam, the Netherlands}
\ead{jellismonaghan@gmail.com}
\ead[url]{https://sites.google.com/site/joellismonaghan/}

\begin{abstract}
Motivated by a problem arising out of DNA origami, we give a general counting framework and enumeration formulas for various cellular embeddings of bouquets and dipoles under different kinds of symmetries. Our algebraic framework can be used constructively to generate desired symmetry classes, and we use Burnside's Lemma with various symmetry groups to derive the enumeration formulas.  Our results assimilate several existing formulas into this unified framework.  Furthermore, we provide new formulas for bouquets with colored edges (and thus for bouquets in nonorientable surfaces) as well as for directed embeddings of directed bouquets.  We also enumerate vertex-labeled dipole embeddings.  Since dipole embeddings may be represented by permutations, the formulas also apply to certain equivalence classes of permutations and permutation matrices. The resulting bouquet and dipole symmetry formulas enumerate structures relevant to a wide variety of areas in addition to DNA origami, including RNA secondary structures, Feynman diagrams, and topological graph theory.
For uncolored objects we catalog $58$ distinct sequences, of which $43$ have not, as far as we know, been described previously.
\end{abstract}


\end{frontmatter}

\part{\Large Overview and results}

\section{Introduction}

We provide a unified framework and new results for enumeration formulas for cellularly embedded bouquets and dipoles under various symmetries.  Bouquets and dipoles are graphs that encode information critical in diverse settings.  A bouquet is a graph with one vertex and some loops.  A dipole is a graph with two vertices and some edges, none of which are loops.  Cellular embeddings of these, or other, graphs are determined by cyclic orderings of half-edges around the vertices, possibly also with edge twists.  We catalog here our own as well as existing enumeration formulas for embedded bouquets and dipoles under various symmetry constraints and equivalences, with the goal of making them readily accessible in one place.  For uncolored objects we list $58$ distinct sequences, of which $12$ already appear in the Online Encyclopedia of Integer Sequences (OEIS) \cite{OEIS}, and three occur elsewhere in the literature.  The remaining $43$ sequences have, as far as we know, not been described previously.

We divide this paper into two parts.  Part I contains some background, definitions, descriptions of the various sorts of embeddings and symmetries, the statements of the counting formulas, and a discussion of open problems.  Part II contains all the technical details of the proofs.  For the convenience of the reader, we list all the formulas in Section \ref{formulas}, deferring the proofs of the results to Sections \ref{pf-dipole}, \ref{pf-bouquet}, and \ref{pf-dirbouquet} in Part II.  We begin with the formulas and proofs for dipoles because the computational ideas there encompass the simpler analogues for bouquets.
Furthermore, we give accessible geometric interpretations and applications for the various symmetries. A reader who simply needs the formulas can go directly to Section \ref{formulas}, although the descriptions of the various objects and symmetries in Section \ref{terminology} may be helpful in identifying the appropriate formulas.

Previous work on enumerating embeddings of graphs, including bouquets and dipoles, under various symmetries has been done by
Mull, Rieper, and White \cite{MRW88},
Rieper \cite{Rie90},
Kwak and Lee \cite{KL94},
Mull \cite{Mul99},
Kim and Park \cite{KP00},
Feng, Kwak, and Zhou \cite{FKZ10, FKZ13}, and
Chen, Gao, and Huang \cite{CGH18}.
We cite results from some of these sources where appropriate in Section \ref{formulas}.
There is also an extensive literature on counting embeddings of bouquets, dipoles and other graphs where symmetry considerations are not taken into account. 
In this situation it is possible to examine the \emph{genus distribution} of embeddings of a given graph, and information on the number of embeddings of a given graph can be represented using a \emph{genus polynomial}.  Early work on genus distributions appears for example in
Stahl \cite{Sta83},
Gross and Furst \cite{GF87},
Furst, Gross, and Statman \cite{FGS89}, 
and Gross, Robbins, and Tucker \cite{GRT89}.
Citation searches on these papers provide access to the recent literature in this area.

\begin{figure}[h]
  \centering
    \includegraphics[clip, trim=2cm 18cm 2cm 3cm, width=0.5\textwidth]%
	{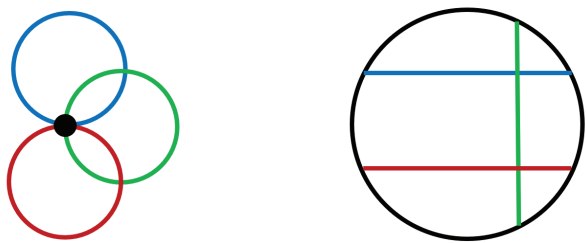}
    \caption{A bouquet (left) and the corresponding chord diagram
        (right).}
    \label{fig:bouquetchord}
\end{figure}

The problem of enumerating bouquets and dipoles arises in surprisingly diverse settings.  For example, bouquets are equivalent to chord diagrams, as in Figure~\ref{fig:bouquetchord}. The cyclic order of the half-edges about the the vertex in the bouquet corresponds to the cyclic order of endpoints of the chords in the chord diagram. Chord diagrams are used in genomics \cite{Circos09} and modeling RNA secondary structures \cite{Bon08, Pen16, Zaj+18}.  Chord diagrams also characterize moduli spaces \cite{AMR96}. Chord diagrams with labeled points and possibly some unpaired points are counted by Feynman integrals. The nLab has an extensive catalog of such applications of chord diagrams in knot theory and physics \cite{nLab}. 

To these many settings for bouquets and dipoles, we add our own motivating application, that of constructing and analyzing DNA origami molecules. Determining routes for a single strand of DNA through assembly targets is integral to both DNA origami~\cite{E+17, E+15cost, V+16} and experimental verification of the targeted constructs~\cite{EE-M19, WJS09}.  When the target  construct has the shape of a graph, these routes correspond to facial walks in an embedding of the graph.
  
In the DNA origami method of self assembly, a single stranded DNA plasmid, called a scaffolding strand,  traces a targeted shape (such a wireframe polyhedron).  Then some 200-250 short strands of DNA complementary to specific regions of the plasmid are introduced to fold and secure the molecule into the desired shape. See \cite{Pel07, Rot06, S15}.

%
\begin{figure}[h]
    \begin{subfigure}{0.2\textwidth}
\centering
    \includegraphics[clip, trim=3cm 12cm 20cm 2.8cm, width=\textwidth]{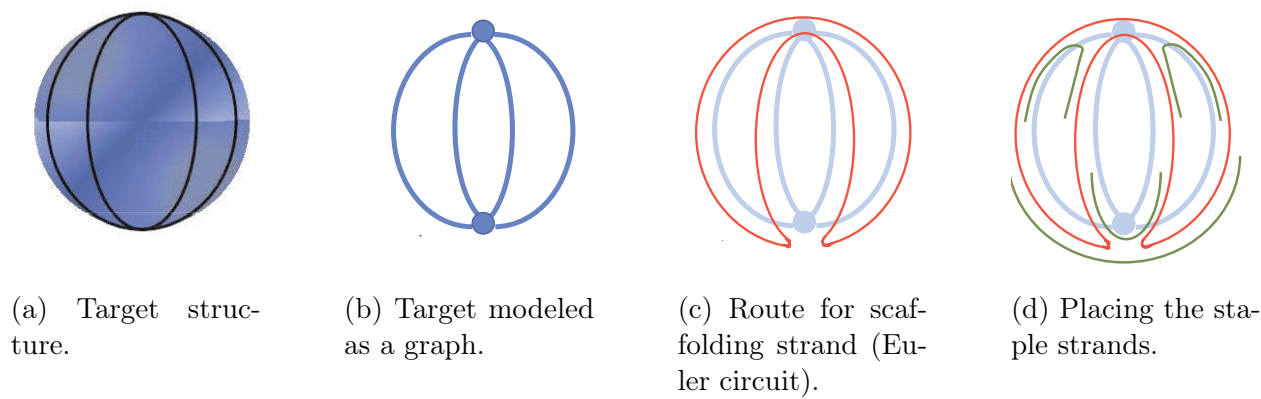}
    \caption{Target structure.\\}
    \label{fig:target}
\end{subfigure}
\hfill
\begin{subfigure}{0.2\textwidth}
\centering
    \includegraphics[clip, trim=9cm 12cm 14cm 2.8cm, width=\textwidth]{figures/TargetScaffold.pdf}
    \caption{Target modeled as a graph.\\}
    \label{fig:targetgraph}
\end{subfigure}
\hfill
\begin{subfigure}{0.2\textwidth}
\centering
    \includegraphics[clip, trim=15cm 12cm 8cm 2.8cm, width=\textwidth]{figures/TargetScaffold.pdf}
    \caption{Route for scaffolding strand (Euler circuit).}
    \label{fig:strand}
\end{subfigure}
\hfill
\begin{subfigure}{0.2\textwidth}
\centering
    \includegraphics[clip, trim=21cm 12cm 2cm 2.8cm, width=\textwidth]{figures/TargetScaffold.pdf}
    \caption{Placing the staple strands.\\}
    \label{fig:staples}
\end{subfigure}
     \caption{Basic design principle for DNA origami.}
    \label{fig:dnadesign}
\end{figure}

When the target is a wireframe structure modeled as as graph embedded in space, a key design step is to determine a route through the graph that the scaffolding strand will follow and to locate the staple strands, as in Figure \ref{fig:dnadesign}.  A related problem, that of finding a reporter strand, corresponds to finding a route through the target graph that traces every edge at least once and at most twice, and when twice in opposite directions.  The resulting route corresponds to a facial walk in a special embedding of a graph, called an edge-outer embedding, as described below.  See~\cite{EE-M19, WJS09}.

\begin{figure}[tb]
\centering
\begin{subfigure}{0.2\textwidth}
\centering
    \vskip15pt
    \includegraphics[clip, trim=2.5cm 20cm 14cm 5cm, width=0.7\textwidth]{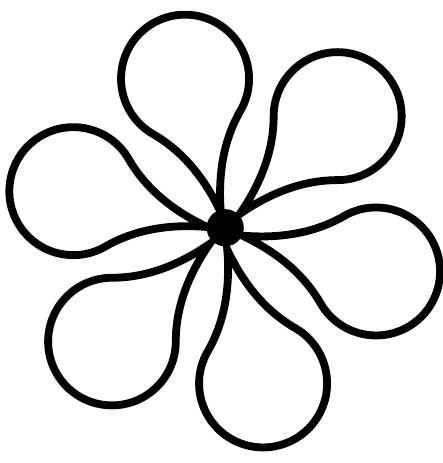}
    \vskip20pt
    \caption{Dual of an upper embeddable  graph with one face.}
    \label{fig:Upper1}
\end{subfigure}
\hfill
\begin{subfigure}{0.2\textwidth}
\centering
    \vskip17pt
    \includegraphics[clip, trim=10cm 19.5cm 4cm 5cm, width=\textwidth]{figures/SurfaceDuals.pdf}
    \vskip20pt
    \caption{Dual of a bi-Eulerian edge-outer embeddable graph.}
    \label{fig:OTEF}
\end{subfigure}
\hfill
\begin{subfigure}{0.2\textwidth}
\centering
    \includegraphics[clip, trim=2.5cm 9cm 14cm 10cm, width=0.5\textwidth]{figures/SurfaceDuals.pdf}
    \caption{Dual of an upper embeddable graph with two faces.}
    \label{fig:Upper2}
\end{subfigure}
\hfill
\begin{subfigure}{0.2\textwidth}
\centering
    \vskip-6pt
    \includegraphics[clip, trim=11.3cm 10cm 3.4cm 10.6cm, width=0.85\textwidth]{figures/SurfaceDuals.pdf}
    \vskip-11pt
    \caption{Dual of an edge-outer embeddable graph.}
    \label{fig:EdgeOuter}
\end{subfigure}
\vskip3pt   
\caption{Surface duals of classes of embedded graphs.}
\label{fig:SurfDual}
\end{figure}

In the theoretical setting of graph embeddings, the bouquets and dipoles characterize some important classes of graphs.  An \emph{upper-embeddable} graph is a graph that can be cellularly embedded in an orientable surface with only one or two faces (see \cite{Xuo79b, Xuo79a}).  An \emph{edge-outer embeddable} graph is a graph that can be cellularly embedded so that every edge lies on a single distinguished face, although there may be other faces as well (see \cite{EE-M19}). These graphs, and their special sub-classes, can be characterized by their surface duals, as shown in Figure \ref{fig:SurfDual}. For simplicity we show planar examples, but in general these embeddings are not necessarily planar.  For example, the half-edges of a loop need not occur consecutively around a vertex.
    
Upper embeddable graphs with exactly one face are characterized by having bouquets as their surface duals, and edge-outer embeddable graphs with exactly two faces that are Euler circuits (\emph{bi-Eulerian} embeddings) are characterized by having dipoles as their surface duals.  To the best of our knowledge there has not yet been any effort to enumerate the looped dipoles in Figure \ref{fig:Upper2} or graphs of the form in Figure \ref{fig:EdgeOuter}.  However, those counting problems are likely to build on the formulas given here, and thus this work lays the necessary foundations for this enumeration.
    
%

In all these settings where bouquets and dipoles play a central role,  good enumeration formulas are essential. In addition to the obvious theoretic interest, these inform experimental design, algorithmic solutions, and estimations of solution space size.

%
\begin{figure}
\centering
\begin{subfigure}{0.4\textwidth}
   \includegraphics[clip, trim=13cm 20cm 1cm 3cm, width=0.75\textwidth]{figures/Linear}
   \caption{A bouquet with fixed labels.}
   \label{fig:labelbouqet}
\end{subfigure} \hfill
\begin{subfigure}{0.4\textwidth}
   \vskip8pt
   \includegraphics[clip, trim=0.5cm 19.5cm 9cm 3cm, width=\textwidth]{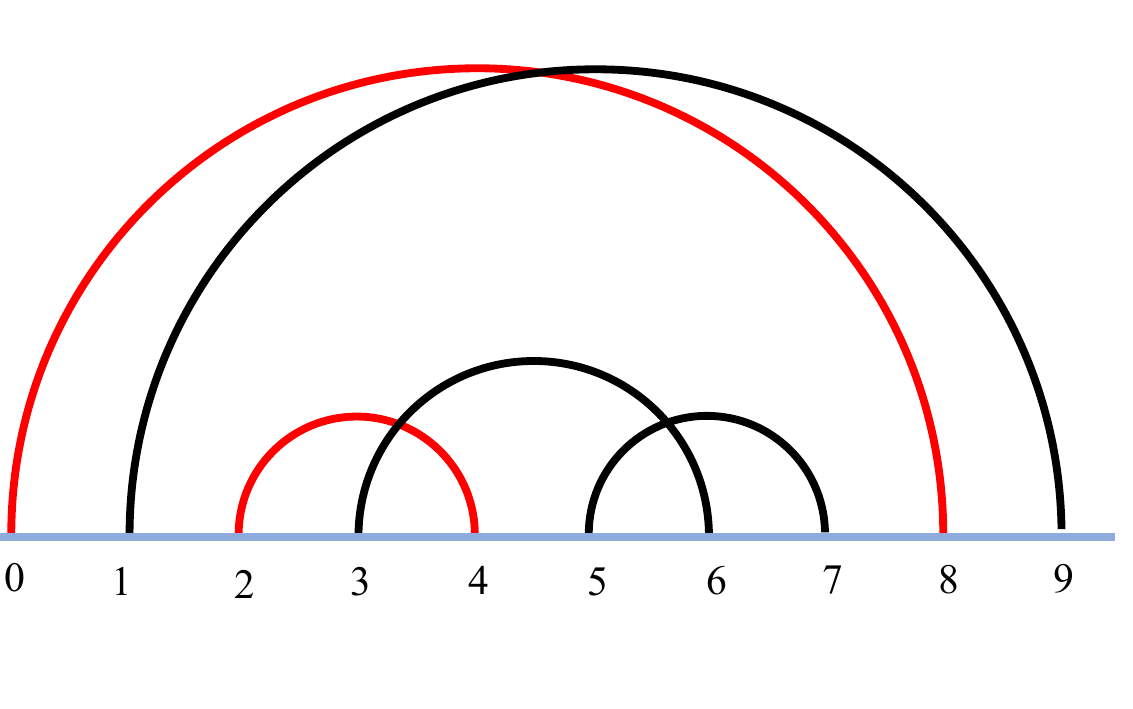}
   \vskip5pt
   \caption{Linear diagram for the bouquet.}
   \label{fig:linebouqet}
\end{subfigure} \hfill
\begin{subfigure}{0.4\textwidth}
   \includegraphics[clip, trim=11cm 6.5cm 4cm 14cm, width=0.65\textwidth]{figures/Linear.pdf}
   \caption{A dipole with fixed labels.}
   \label{fig:labeldipole}
\end{subfigure} \hfill
\begin{subfigure}{0.4\textwidth}
   \includegraphics[clip, trim=1cm 6cm 12cm 14cm, width=0.8\textwidth]{figures/Linear.pdf}
   \caption{Linear diagram for the dipole.}
   \label{fig:LineDipole}
\end{subfigure}
\vskip3pt
\caption{Bouquets and dipoles with fixed labels
have no rotational symmetry, so can be identified with linear diagrams.}
\label{fig:Linearized}
\end{figure}

However, each setting requires careful consideration of what symmetries are relevant to the application.  For example, linear RNA secondary structures are often modeled by chord diagrams with a designated point on the circle boundary that indicates where the chord diagram should be `cut open' to form a linear structure. See Figure \ref{fig:Linearized}.
The presence or absence of such a symmetry-breaking point significantly alters the counting problem. Similarly, chirality often plays a role, and this too changes the enumeration problem.  The possibility of directions on the edges, or of edges with various attributes (here colors), lead to further enumeration problems.

Because of the diverse applications, existing formulas are widely scattered in the literature, making them sometimes challenging to find. Thus, we include known results here as well, incorporating them into our overall framework and in some cases simplifying the proofs or formulas. We then complete the work of finding enumeration formulas for the remaining symmetries in the groups we consider, particularly providing enumeration formulas for bouquets with colored or directed edges and for orientable (as opposed to oriented) embeddings of dipoles.  When the number of colors is two, the formulas for colored bouquets enumerate embeddings of bouquets in nonorientable surfaces.  Greater numbers of colors can be used in applications to differentiate types of edges, for example to differentiate different types of attachments in RNA secondary structures.  The enumeration of dipolar cogs and other structures related to dipole embeddings also uses some ideas that may support the more challenging problem of enumerating dipole embeddings in nonorientable surfaces.

We have verified our counting formulas for small values of $n$ (the number of edges) and, where appropriate, $k$ (the number of colors) by explicit construction by computer of the objects being counted.
Specifically, we verified formulas (D1)--(D8), (B1)--(B2), and (A1)--(A5) in Section \ref{formulas}.  Since all other formulas here are linear combinations of these, this also provides verification for our other results.
Our programs use the same framework as the counting formulas.  We consider our objects as orbits of easily-generated basic objects under various group actions. 
While generating the basic objects in a fixed order, the programs test whether a given object is the earliest in its orbit under action of the group elements, thereby identifying a unique representative of each orbit.

\section{Terminology and notation}\label{terminology}

\subsection{Embedding concepts}\label{term-obj}

The objects we are counting are embeddings of graphs in compact surfaces, or objects related to these. We assume that the reader is familiar with embeddings of graphs and their combinatorial representations; standard references are \cite{GT, MT}, and details for the \emph{cogs} defined below may be found in \cite{E-MM13} .
Unless stated otherwise, all embeddings of graphs and digraphs in this paper will be \emph{cellular}, meaning that each face is homeomorphic to an open disk.

An embedding of a graph in an orientable surface can be described up to homeomorphism by giving a \emph{rotation scheme} specifying a \emph{rotation}, that is, a cyclic ordering of the half-edges, at each vertex.  A \emph{generic} embedding (in either an orientable or nonorientable surface) of a graph can be described up to homeomorphism by a rotation scheme together with \emph{edge signatures} specifying whether each edge is \emph{twisted} or \emph{untwisted}.  The representation using a rotation scheme is unique for oriented embeddings but not for orientable embeddings (see below for the distinction).  The representation of a generic embedding using a rotation scheme and edge signatures is in general not unique.
 
An \emph{oriented embedding} of a graph is an embedding in an orientable surface with a definite clockwise orientation, up to graph isomorphism and orientation-preserving surface homeomorphism.
Oriented embeddings are in one-to-one correspondence with rotation schemes (without edge signatures).
The \emph{reflection} of an oriented embedding is obtained by reversing the clockwise orientation of the surface.
An \emph{orientable embedding} of a graph is an embedding in an orientable surface where the clockwise orientation is not specified, up to graph isomorphism and surface homeomorphism.  Equivalently, orientable embeddings are equivalence classes of oriented embeddings under reflection.

We make an explicit distinction between oriented and orientable embeddings.  In many settings the distinction between a surface being oriented and being orientable is either unnecessary or implicitly understood.  However, for enumerative results it is important to distinguish between these.

We also consider \emph{cogs}, also known as \emph{cyclically ordered graphs} or \emph{rigid-vertex graphs}.  These are graphs with a rotation at each vertex as above, but here two cogs are equivalent if there is a graph isomorphism between the two underlying graphs which at each vertex either preserves the rotation at the vertex or reverses it.  In other words, there is an undirected cyclic ordering, defined only up to reversal, at each vertex, instead of the directed cyclic ordering in a rotation scheme.
The edges of cogs have no signatures.
Cogs are important because a cog represents an equivalence class of graph embeddings under partial Petrie duality (edge twisting) operations (see \cite[Lemma 3.16(2)]{E-MM13}), or an equivalence class of orientable embeddings under vertex flips (rotation reversals).

A \emph{digraph} is a graph with a direction (from one end-vertex to the other) specified for each edge.  The directed edges are called \emph{arcs}.  A \emph{directed embedding} of a digraph is an embedding where every facial walk is a directed walk.  This is equivalent to the property that at each vertex the half-arcs in the rotation alternate in direction between outwards and inwards.

Objects based on graphs or digraphs can be considered to be \emph{vertex-labeled}, if each vertex has a unique label that must be preserved by any symmetry operation (although edges can be permuted), or \emph{vertex-unlabeled} if symmetry operations that permute vertices are allowed.
In this paper we deal with bouquets, where this distinction is irrelevant, and with dipoles, where we have only two vertices, which can be swapped in the vertex-unlabeled situation.  We will say explicitly if an object derived from a dipole is vertex-labeled; otherwise it is assumed to be vertex-unlabeled.

\subsection{Bouquets and dipoles}\label{bouquetdipole}

A \emph{bouquet} is a graph with exactly one vertex, and a \emph{dipole} is a graph with exactly two vertices and no loops.  These may be embedded in either orientable or nonorientable surfaces, or given related structures, such as a cog structure. We will consistently denote the number of edges by $n$.

In our work we will label the half-edges around each vertex in an embedded bouquet or dipole to encode embedding information.  Enumerating embeddings with fixed half-edge labels (see the discussion of `(colored) labeled bouquets' and `labeled dipoles' in Part II) is easy, since only the identity permutation leaves the fixed labels in their original position.  Since fixed edge labels break cyclic symmetries, the diagrams maybe `linearized' at the vertices, that is, starting at the half-edges labeled $0$ and opening up the vertices to lines, as in Figure \ref{fig:Linearized}. The number of labeled bouquets is the same as the number of chord diagrams, which is the number of perfect matchings in a complete graph of order $2n$, namely $(2n-1)!! = (2n-1)(2n-3)(2n-5) \dots \cdot 3 \cdot 1$.  The number of labeled dipoles is just $n!$.

Since bouquets and dipoles with fixed labels have no symmetries and their counting formulas are easy and well known, we focus our attention on bouquets and dipoles where the labels may be permuted.  Here, the permutations reveal the underlying structural symmetries of the embedded graph, and the enumeration formulas are much more complex. 

\begin{figure}
\centering
\begin{subfigure}{0.44\textwidth}
\centering
    \null\hskip0.02\textwidth
    \includegraphics[clip, trim=12cm 18cm 0cm 2.5cm, width=\textwidth]{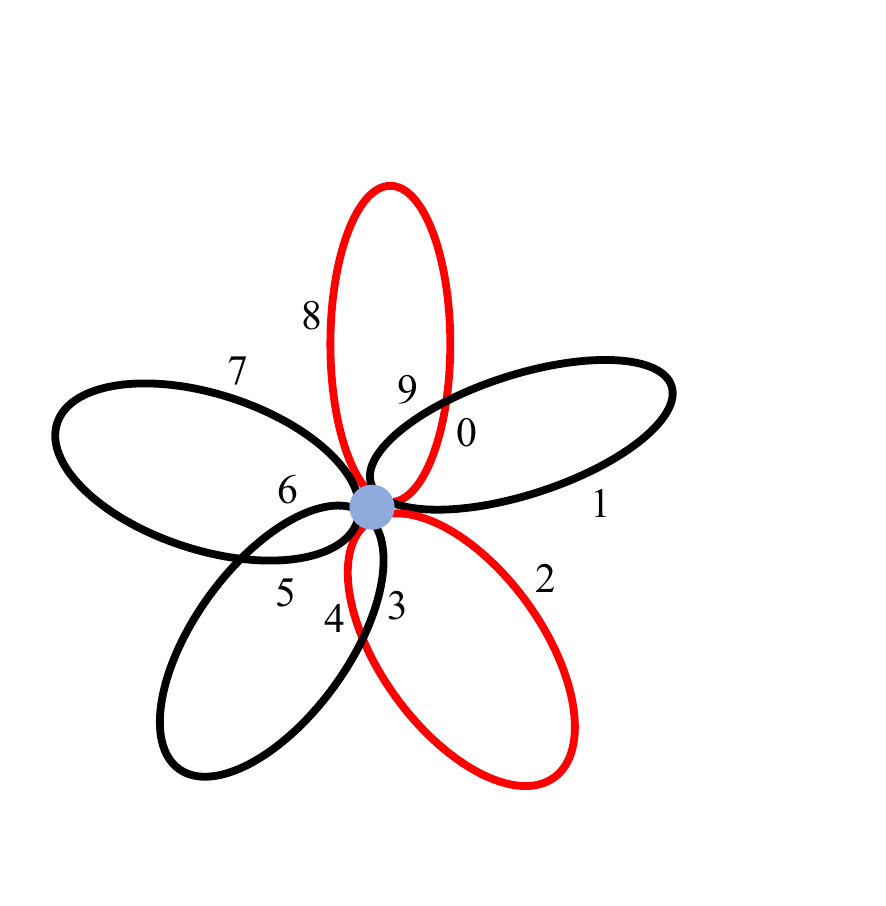}
    \null\hskip0.02\textwidth
    \caption{A colored bouquet embedding with half-edges labeled in
	rotational order.}
    \label{fig:InterleaveII}
\end{subfigure}
\break
\begin{subfigure}{0.44\textwidth}
    \null\hskip0.02\textwidth
    \includegraphics[clip, trim=1cm 18cm 11.75cm 3.2cm, width=0.8\textwidth]{figures/Interleaving.pdf}
    \null\hskip0.02\textwidth
    \caption{The colored labeled abstract bouquet.}
    \label{fig:InterleaveI}
\end{subfigure}
\hfill
\begin{subfigure}{0.44\textwidth}
    \null\hskip0.04\textwidth
    \includegraphics[clip, trim=8.25cm 9cm 6cm 13cm, width=0.7\textwidth]{figures/Interleaving.pdf}
    \null\hskip0.0\textwidth
    \vskip9pt
    \caption{The colored chord diagram.}
    \label{fig:InterleaveChord}
\end{subfigure}  
\vskip3pt
\caption{Colored labeled abstract bouquets or chord diagrams correspond to
colored embedded bouquets labeled in rotational order.}
\label{fig:Interleaving}
\end{figure}

We can convert an embedding of a bouquet (possibly with colors or edge directions) into a labeled bouquet, an abstract graph with labels on the half-edges.
We label the half-edges of the embedding with elements of the cyclic group $\Ztn$ in order around the vertex.  
The two labels on each edge encode the position of that edge in the embedding. We can then simply consider the underlying abstract graph, with these half-edge labels.
For example, the colored bouquet embedding represented by the rotation shown in Figure \ref{fig:InterleaveII} is represented by the colored labeled bouquet in Figure \ref{fig:InterleaveI} (drawn in the plane for convenience, although only the graph structure and labels matter).
We can also draw a chord diagram, which can be thought of as an `inside-out' drawing of the embedded bouquet.  We expand the vertex of the bouquet into a large circle, and draw the edges inside this circle, as chords, instead of outside it.  We transfer the label on each half-edge to the point at which that half-edge meets the circle, so that the points are labeled by elements of $\Ztn$ in cyclic order around the circle.  The chords form a \emph{perfect matching}, a partition of all of the vertices into pairs, in the complete graph whose vertices are the elements of $\Ztn$.
Figure \ref{fig:InterleaveChord} shows the chord diagram corresponding to the other two parts of Figure \ref{fig:Interleaving}.

%
We enumerate embeddings of bouquets up to the following symmetries.
We consider embeddings of bouquets in an oriented surface up to rotational symmetry.  Thus, two bouquets that can be superimposed by rotating one of them are considered the same, but reflection is not allowed.
Allowing reflections (reversal of the orientation of the surface) means we are counting embeddings in an orientable surface, where we allow orientation-reversing homeomorphisms of the surface.  We also consider the reflexible embeddings, preserved by reflection, and the chiral embeddings, which are not preserved by reflection.

One of our main contributions for bouquets is counting embeddings of colored bouquets, which also allows us to count generic (orientable and nonorientable) embeddings of bouquets, and hence, by simply subtracting, the nonorientable embeddings of bouquets.  We also count directed embeddings of \emph{directed bouquets}, i.e., directed graphs with one vertex.  In this setting we again count embeddings with colors, which allows us to count generic and nonorientable embeddings.

%
\begin{figure}[h]
\centering
\begin{subfigure}{0.375\textwidth}
\centering
    \includegraphics[clip, trim=10cm 18cm 2cm 2.5cm, width=\textwidth]{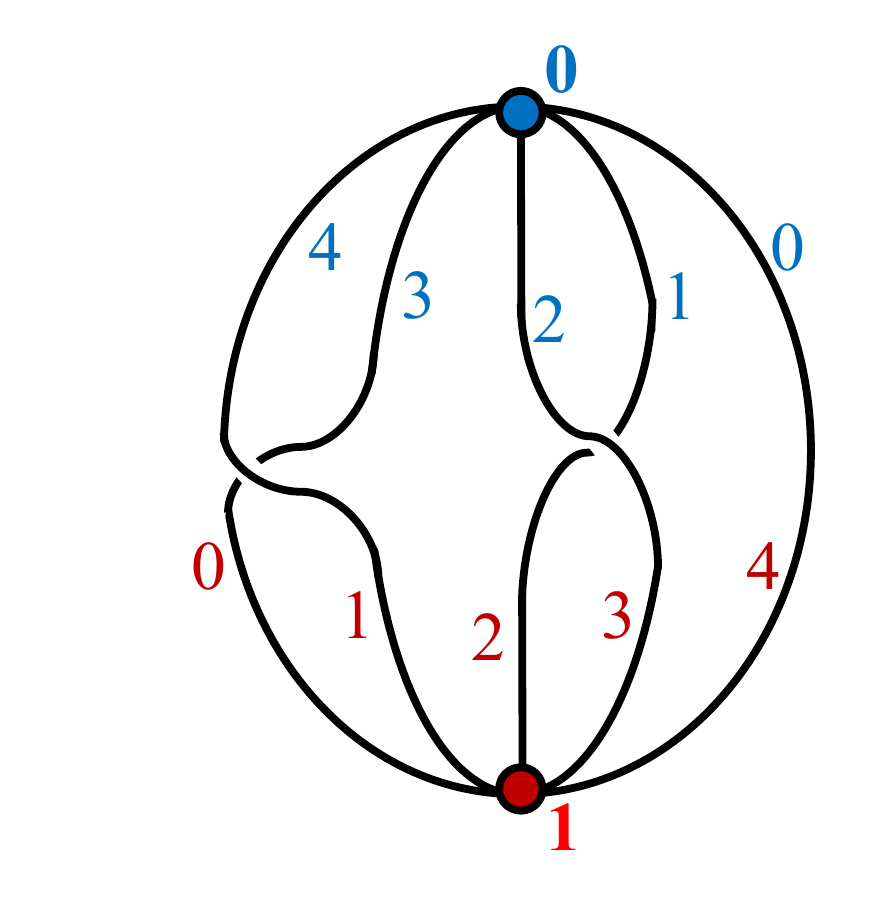}
    \caption{Dipole embedding with half-edges labeled in rotational
            order.}
    \label{fig:dipolesI}
\end{subfigure}
\hfill
\begin{subfigure}{0.4\textwidth}
    \includegraphics[clip, trim=1cm 17.3cm 10.5cm 2.7cm, width=\textwidth]{figures/dipoles.pdf}
    \caption{Corresponding labeled dipole.}
    \label{fig:dipolesII}
\end{subfigure}
\vskip 3pt
\caption{Labeled dipoles correspond to a dipole embedding.}
\label{fig:dipoles}
\end{figure}

We enumerate embeddings of dipoles up to the following symmetries. 
Since dipoles have two vertices, we can consider the vertices as being distinct, e.g., having two different labels,  or not. 
For embeddings in oriented surfaces we can rotate (cyclically shift) the labels around each vertex independently.
For embeddings in orientable surfaces we also allow reflection, i.e., reversal of the orientation of the surface.
For dipoles, we also count cogs, both with and without distinguished vertices.

Our approach for counting embeddings of dipoles will be similar to that for counting embeddings of bouquets.  We turn a dipole embedding into an object called a labeled dipole, by labeling the edges in rotational order around each vertex.  See Figure \ref{fig:dipoles}, where the rotational information for the embedding in Figure \ref{fig:dipolesI} is represented by the correspondence between labels around the two vertices in Figure \ref{fig:dipolesII}.  See also Figure \ref{fig:dipoleperms} later, which shows how an embedding may be represented by more than one labeled dipole.

\begin{figure}
\centering
\begin{subfigure}{0.25\textwidth}
    \includegraphics[clip, trim=0.5cm 19.5cm 15cm 4cm, width=\textwidth]{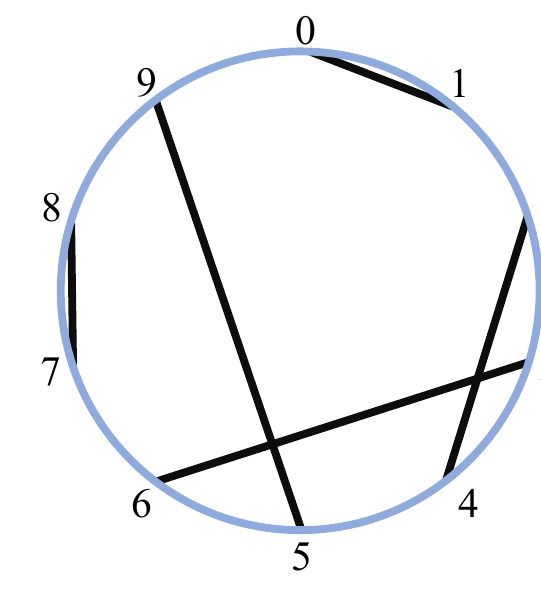}
    \caption{Neither reflexive nor rotational symmetry.\\}
    \label{fig:BSI}
\end{subfigure}
\hfill
\begin{subfigure}{0.25\textwidth}
\centering
    \includegraphics[clip, trim=8cm 19.5cm 7.5cm 4cm, width=\textwidth]{figures/bouquetsymmetries.pdf}
    \caption{$180^{\circ}$ rotational symmetry; no reflexive symmetry.}
    \label{fig:BSII}
\end{subfigure}
\hfill
\begin{subfigure}{0.25\textwidth}
\centering
    \includegraphics[clip, trim=15cm 19.5cm 0.5cm 4cm, width=\textwidth]{figures/bouquetsymmetries.pdf}
    \caption{Reflexive symmetry through 0-5; no rotational symmetry.}
    \label{fig:BSIII}
\end{subfigure}
\break
\vskip5pt
\begin{subfigure}{0.25\textwidth}
    \includegraphics[clip, trim=1cm 12cm 14.5cm 12cm, width=\textwidth]{figures/bouquetsymmetries.pdf}
    \caption{$180^{\circ}$ rotation plus arc-reversal; reflection through 4-9; reflection through $1\frac12 - 6\frac12$ plus arc-reversal.}
    \label{fig:BSIV}
\end{subfigure}
\hfill
\begin{subfigure}{0.25\textwidth}
\centering
    \includegraphics[clip, trim=15cm 12cm 0.5cm 12cm, width=\textwidth]{figures/bouquetsymmetries.pdf}
    \caption{Rotation by $72^{\circ}$; reflection through $\frac12 -
5\frac12$ plus arc-reversal.\\ \\}
    \label{fig:BSV} 
    \end{subfigure}
    \hfill
\begin{subfigure}{0.25\textwidth}
\centering
    \includegraphics[clip, trim=8.5cm 12cm 7cm 12cm, width=\textwidth]{figures/bouquetsymmetries.pdf}
    \caption{Neither reflexive nor rotational symmetry, although the
underlying uncolored diagram has many symmetries.}
    \label{fig:BSVI}
\end{subfigure}
\vskip3pt
\caption{Chord diagram symmetries.}
\label{fig:chordsymm}
\end{figure}

The difficulty in counting comes from the fact that different embeddings correspond to different numbers of basic labeled objects (for us, colored labeled bouquets or labeled dipoles).  
To use embeddings of bouquets as an example, the number of labeled bouquets or chord diagrams that correspond to a given bouquet embedding depends on the intrinsic symmetries of the chord diagrams, and also on which symmetries matter for the kind of embedding we are counting.  Figure \ref{fig:chordsymm} shows a number of different symmetries.  If we have an oriented embedding, where we only allow rotational symmetry, all rotations of Figure \ref{fig:BSI} give different chord diagrams, so there are $10$ different chord diagrams for the corresponding embedding, but rotations of Figure \ref{fig:BSII} by $5$ (corresponding to $180^\circ$) give the same chord diagram, so there are only $5$ different chord diagrams for the corresponding embedding. If we are counting orientable embeddings rather than oriented embeddings, we also need to consider reflexive symmetry, as in Figure \ref{fig:BSIII}.
The other parts of Figure \ref{fig:chordsymm} show examples of other symmetries that can affect counting for different types of objects.  The tool we use to handle this difficulty is Burnside's Lemma, Lemma \ref{burnside}.

%

\subsection{Symmetry and asymmetry under an involution}\label{term-inv}

In many of our results we count a set of objects, and we have an involution (i.e., self-inverse permutation) that can be applied to these objects.  If we count the number of equivalence classes of our objects under the involution, we can also easily count the number of symmetric objects (preserved by the involution), and the number of pairs of the remaining asymmetric objects, where the objects in each pair are swapped by the involution.

A common type of involution is some kind of reflection, and we would like to count the \emph{reflexible} or \emph{achiral} objects that are preserved by reflection, and the number of pairs of \emph{chiral} objects that are asymmetric under (i.e., not preserved by) reflection.  To give a specific example, oriented embeddings of a graph can be divided into reflexible embeddings and chiral pairs of embeddings, where the reflection operation is the reversal of the clockwise orientation of the surface.

The following general formulas apply for a given involution.
\begin{align*}
\text{\#objects} &= \text{\#symmetric} + 2\text{\#(asymmetric pairs)},
		\\
\text{\#(involution classes)} &= \text{\#symmetric} +
\text{\#(asymmetric pairs)}.
\end{align*}
These equations mean that any two of these four numbers determine the other two.  In our results, generally we can count the number of objects and the number of involution classes of objects; we can then compute the number of symmetric objects and pairs of asymmetric objects as follows.
\begin{align}
\text{\#symmetric} &= 2\text{\#(involution classes)} -
\text{\#objects},
	\label{symmetric}\\
\text{\#(asymmetric pairs)} &= \text{\#objects} - \text{\#(involution
	classes)}.
	\label{asymmetric}
\end{align}
We will apply these formulas many times.

Besides reflection, we have other types of involutory symmetry,
If our objects are vertex-labeled dipoles and our involution is exchanging the two vertices, we refer to the equivalence classes as \emph{vertex-unlabeled dipoles} or just \emph{dipoles}. The symmetric ones are \emph{vertex-interchangeable} and the asymmetric ones come in \emph{non-vertex-interchangeable} pairs.
If our objects are directed graphs and the involution is reversing all of the arcs, we refer to the equivalence classes as \emph{arc-reversal classes}. The symmetric ones are \emph{arc-reversible} and the asymmetric ones come in \emph{arc-irreversible} pairs.

Although it should be intuitively clear that our specific involutions make sense, we provide some technical details for the interested reader. Generally we have a set of basic objects $\cS$ and a group $\Ga$ acting on $\cS$.  If we have a normal subgroup $\De \nsubeq \Ga$, then elements of $\Ga$ act on $\cS/\De$, the orbits of $\cS$ under the action of $\De$, which we consider as a set of derived objects.  In particular, if $\De$ is a subgroup of $\Ga$ and $|\Ga| = 2|\De|$ then $\De \nsubeq \Ga$, and if $J \in \Ga-\De$ is an involution in $\Ga$ then $J$ acts as an involution on the set of derived objects $\cS/\De$.
When we use the results of this subsection the situation is always an instance of this general framework.

\subsection{Standard counting functions}\label{term-func}

In our results we often use Euler's totient function $\eul(n)$, which is the number of integers $k$ that are relatively prime to $n$ and satisfy $1 \le k \le n$.
We use the standard notation $(a,b)$ for the greatest common divisor of $a$ and $b$.

We also frequently use $\ma(n, j)$, the number of $j$-matchings (sets of $j$ disjoint $2$-subsets) of an $n$-set, which is
$$\ma(n,j) = \binom{n}{2j} \frac{(2j)!}{ 2^j j!}
	= \binom{n}{2j} (2j-1)!!
	= \frac{n!}{(n-2j)!\, 2^j j!}
$$
where $(2j-1)!! = (2j-1)(2j-3)(2j-5) \dots 1$ is a double factorial.
To keep our formulas simple, we often keep $\ma$ in our final expressions rather than replacing it with expressions involving factorials and powers.

\section{The enumeration formulas}\label{formulas}

In this section we summarize our results.  For results that are already known we provide references, and in particular references to entries in the Online Encyclopedia of Integer Sequences \cite{OEIS}.  If a result is not attributed, it is (as far as we know) new, although some results can be derived from known results using the approach in Subsection \ref{term-inv}.

We begin with dipole results, then present the (undirected) bouquet results, and finally the directed bouquet results. In each section we give a small number of basic quantities (technically \emph{coset averages}, defined in Part II) and then all of the counting results are expressed as simple linear combinations of these quantities.   The formulas are indexed by (D1), (D2), $\ldots$ for the dipoles, (B1), (B2), $\ldots$ for the bouquets, and (A1), (A2), $\ldots$ for directed bouquets.  These indices correspond to the proofs of the formulas given in Sections \ref{bo-pf-dipole}, \ref{bo-pf-bouquet}, and  \ref{bo-pf-dirbouquet}, respectively.    For each enumeration formula, we  we give the sequence of values for $n$ (the number of edges) with $0 \le n \le 12$.  For colored objects we just provide the values for $k$ (number of colors) equal to $1$, i.e., in the uncolored situation.

The formulas we give below are generally valid only for $n \ge 1$. For $n=0$ there is always one trivial object with no edges, which is preserved under all symmetries.  Therefore, the numbers for $n=0$ are either $1$ (for objects under some equivalence relation, or objects symmetric under some involution) or $0$ (for pairs of objects asymmetric under an involution).

\subsection{Dipole results}\label{res-dipole}

\subsubsection{Dipole coset averages}\label{ca-res-dipole}

There are five basic quantities for counting dipole embeddings and related objects.  The superscripts on $\delta$ are mnemonics for the reflections, rotations, and vertex exchanges, as detailed in Part II.
As we show in Part II, all the following functions have integer values, given positive integer inputs.
For $n=0$ the value of all of these should be taken to be $1$; the following formulas apply for $n \ge 1$.
\begin{align*}
 \Adx(I)
   &= \frac{1}{n} \sum_{(d,g)\,:\, dg = n} \eul(d)^2\, (g-1)!\,
d^{g-1};
\displaybreak[0]
\\[4pt]
\Adx(X)
   &= \frac{1}{n} \sum_{\substack{(d,g)\,:\, dg = n \\ \text{$d, g$
even}}}
		\eul(d)\, \ma(g, g/2)\, d^{g/2}
	+ \frac{1}{n} \sum_{\substack{(d,g)\,:\, dg = n \\ \text{$d$ odd}}}
		\eul(d) \sum_{j=0}^{\lfloor g/2 \rfloor} \ma(g, j)\, d^j \\*
    & \qquad\text{(the first term here is nonzero only if $n \equiv 0$
	\ (mod $4$))};
\displaybreak[0]
\\[4pt]
 \Adx(R) &= \begin{cases}
	\displaystyle \left(\frac{n-1}{2}\right)!\; 2^{(n-1)/2}
		& \text{for odd $n$}, \\
 	\displaystyle (n+2) \left(\frac{n}{2}-1\right)!\; 2^{n/2-3}
		& \text{for even $n$};
  \end{cases}
\displaybreak[0]
\\[4pt]
 \Adx(R_1)
   &= \begin{cases}
     1 & \text{if $n=1$ or $2$,}\\
     0 & \text{if $n \ge 3$ is odd,}\\
 \displaystyle \frac{1}{n} \left( \frac{n}{2} \right)! \; 2^{n/2-1}
	& \text{if $n\ge 4$ is even;}
   \end{cases}
\displaybreak[0]
\\[4pt]
 \Adx(R_1 X)
   &= \begin{cases}
	\ma(n/2, n/4)\, 2^{n/4-1}
		& \text{if $n \equiv 0 \pmod4$,}\\
	\ma((n-1)/2, (n-1)/4)\, 2^{(n-1)/4}
		& \text{if $n \equiv 1 \pmod4$,}\\
	\ma((n-2)/2, (n-2)/4)\, 2^{(n-2)/4}
		& \text{if $n \equiv 2 \pmod4$,}\\
	0
		& \text{if $n \equiv 3 \pmod4$.}
    \end{cases}
\end{align*}
The values of $\Adx(I)$ and $\Adx(R_1)$ occur in the OEIS \cite{OEIS}; see (D1) and (D3/D5:S) below.

\subsubsection{Counting basic dipole objects}\label{bo-res-dipole}

We note that dipoles can be represented using permutations, as we discuss in more detail in Section \ref{pf-dipole} of Part II.  Therefore, most of the results that we mention have natural interpretations in terms of permutations, or equivalently \emph{permutation matrices}, with all entries $0$ except for one $1$ in each row and each column.

For permutations, specifically elements of the symmetric group $\Sz_n$, the operations defined in Section \ref{pf-dipole} act by  cyclically shifting the input variable ($S_0$), cyclically shifting the output variable ($S_1$), reversing the input variable ($R_0$), reversing the output variable ($R_1$), simultaneously reversing the input and output variables ($R$), and inverting the permutation ($X$).  For formulas (D1)--(D6) we briefly indicate what they count in terms of equivalence classes of permutations under these operations.

\begin{countitem}
(D1) The number of vertex-labeled oriented $n$-edge dipole embeddings is $\Adx(I)$.

This also counts elements of $\Sz_n$ equivalent under $S_0$ and $S_1$.
It appears in the OEIS \cite{OEIS} as A002619.

\values D1 1,  1, 1, 2,  3, 8, 24,  108, 640, 4492,  36336, 329900, 3326788.
\end{countitem}

\begin{countitem}
(D2) The number of oriented $n$-edge dipole embeddings is
$\rc2 (\Adx(I)+\Adx(X))$.

This also counts elements of $\Sz_n$ equivalent under $S_0$, $S_1$ and $X$.
This sequence was found by Rieper \cite[Theorem 5.10]{Rie90} and also by Feng, Kwak, and Zhou \cite[Theorem 4.1]{FKZ10}.

\values D2 1,  1, 1, 2,  3, 7, 19,  71, 369, 2393,  18644, 166573, 1669243.
\end{countitem}

\begin{countitem}
(D3) The number of vertex-labeled orientable $n$-edge dipole embeddings is
$\rc2 (\Adx(I)+\Adx(R))$.

This also counts elements of $\Sz_n$ equivalent under $S_0$, $S_1$, and $R$.

\values D3  1,  1, 1, 2,  3, 8, 20,  78, 380, 2438,  18744, 166870, 1670114.
\end{countitem}

\begin{countitem}
(D4) The number of orientable $n$-edge dipole embeddings is
$\rc4 (\Adx(I) + \Adx(R) + 2\Adx(X))$.

This also counts elements of $\Sz_n$ equivalent under $S_0$, $S_1$, $R$, and $X$.

\values D4   1,  1, 1, 2,  3, 7, 17,  56, 239, 1366,  9848, 85058, 840906.
\end{countitem}

\begin{countitem}
(D5) The number of vertex-labeled $n$-edge dipolar cogs is
$\rc4 (\Adx(I) + \Adx(R) + 2\Adx(R_1))$.

This also counts elements of $\Sz_n$ equivalent under $S_0$, $S_1$, $R_0$, and $R_1$.
It appears in the OEIS \cite{OEIS} as A000940 (with some initial terms missing).

\values D5   1,  1, 1, 1,  2, 4, 12,  39, 202, 1219,  9468, 83435, 836017.
\end{countitem}

\begin{countitem}
(D6) The number of $n$-edge dipolar cogs is
$\rc8 (\Adx(I) + \Adx(R) + 2\Adx(R_1) + 2\Adx(X) + 2\Adx(R_1X))$.

This also counts elements of $\Sz_n$ equivalent under $S_0$, $S_1$, $R_0$, $R_1$, and $X$.
It appears in the OEIS \cite{OEIS} as A006841.

\values D6   1,  1, 1, 1,  2, 4, 10,  28, 127, 686,  4975, 42529, 420948.
\end{countitem}

Our final items in this subsubsection do not have natural interpretations in terms of dipole embeddings, but can be expressed in terms of permutations or permutation matrices.

\begin{countitem}
(D7) The number of equivalence classes of permutations in $\Sz_n$ under cyclic shifts of input variable ($S_0$), cyclic shifts of output variable ($S_1$), and reversal of input variable only ($R_0$) is 
$\rc2 (\Adx(I)+\Adx(R_1))$.

This also counts elements of $\Sz_n$ equivalent under $S_0$, $S_1$, and $R_1$.  
It appears in the OEIS \cite{OEIS} as A000939.

\values D7   1,  1, 1, 1,  2, 4, 14,  54, 332, 2246,  18264, 164950, 1664354.
\end{countitem}

\begin{countitem}
(D8) The number of equivalence classes of $\nn$ permutation matrices under cyclic shifts of the row set, cyclic shifts of the column set, and rotation of the matrix by multiples of $90^\circ$ is
$\rc4( \Adx(I)+\Adx(R)+2\Adx(R_1 X) )$.

\values D8  1,  1, 1, 1,  2, 5, 11,  39, 193, 1225,  9378, 83435, 835087.
\end{countitem}

\subsubsection{Symmetric and asymmetric dipole
objects}\label{sa-res-dipole}

Here we provide some formulas obtained by applying equations \eqref{symmetric} and \eqref{asymmetric} of Subsection \ref{term-inv}.
Given an involution that acts on objects counted by (D$i$) and creates equivalence classes counted by (D$j$), item (D$i$/D$j$) states the formulas for the numbers of symmetric objects, (S) or (D$i$/D$j$:S), and asymmetric pairs of objects, (AP) or (D$i$/D$j$:AP).

First we consider situations where our involution is `reflection' in the sense of reversal of the surface orientation.

\begin{countitem}(D1/D3) For vertex-labeled oriented $n$-edge dipole embeddings
the number of reflexible ones (S) is
$2\text{(D3)}-\text{(D1)} = \Adx(R)$, and
the number of chiral pairs (AP) is
$\text{(D1)}-\text{(D3)} = \rc2 (\Adx(I)-\Adx(R))$.

\svalues S:D1/D3   1,  1, 1, 2,  3, 8, 16,  48, 120, 384,  1152, 3840, 13440.
\apvalues AP:D1/D3  0,  0, 0, 0,  0, 0, 4,  30, 260, 2054,  17592, 163030, 1656674.
\end{countitem}

\begin{countitem}(D2/D4) For oriented $n$-edge dipole embeddings
the number of reflexible ones (S) is
$2\text{(D4)}-\text{(D2)} = \rc2 (\Adx(R)+\Adx(X))$, and
the number of chiral pairs (AP) is
$\text{(D2)}-\text{(D4)} = \rc4 (\Adx(I)-\Adx(R))$ (which is half of (D1/D3:AP)).

The number of reflexible ones (S) was found by Feng, Kwak, and Zhou \cite[Theorem 5.1]{FKZ13}.

\svalues S:D2/D4   1,  1, 1, 2,  3, 7, 15,  41, 109, 339,  1052, 3543, 12569.
\apvalues AP:D2/D4  0,  0, 0, 0,  0, 0, 2,  15, 130, 1027,  8796, 81515, 828337.
\end{countitem}

Next we consider situations where our involution is exchanging vertices.

\begin{countitem}
(D1/D2) For oriented $n$-edge dipole embeddings
the number of vertex-interchangeable ones (S) is
$2\text{(D2)}-\text{(D1)} = \Adx(X)$, and
the number of non-vertex-interchangeable pairs (AP) is
$\text{(D1)}-\text{(D2)} = \rc2 (\Adx(I)-\Adx(X))$.

\svalues S:D1/D2   1,  1, 1, 2,  3, 6, 14,  34, 98, 294,  952, 3246, 11698.
\apvalues AP:D1/D2  0,  0, 0, 0,  0, 1, 5,  37, 271, 2099,  17692, 163327, 1657545.
\end{countitem}

\begin{countitem}
(D3/D4) For orientable $n$-edge dipole embeddings
the number of vertex-interchangeable ones (S) is
$2\text{(D4)}-\text{(D3)} = \Adx(X)$ (equal to (D1/D2:S), see values above), and
the number of non-vertex-interchangeable pairs (AP) is
$\text{(D3)}-\text{(D4)} = \rc4 (\Adx(I)+\Adx(R)-2\Adx(X))$.

\smallskip
\apvalues AP:D3/D4  0,  0, 0, 0,  0, 1, 3,  22, 141, 1072,  8896, 81812, 829208.
\end{countitem}

\begin{countitem}
(D5/D6) For $n$-edge dipolar cogs
the number of vertex-interchangeable ones (S) is
$2\text{(D6)}-\text{(D5)} = \rc2 (\Adx(X)+\Adx(R_1X))$, and
the number of non-vertex-interchangeable pairs (AP) is
$\text{(D5)}-\text{(D6)} =
	\rc8 (\Adx(I)+\Adx(R)+2\Adx(R_1)-2\Adx(X)-2\Adx(R_1X))$.

\svalues S:D5/D6   1,  1, 1, 1,  2, 4, 8,  17, 52, 153,  482, 1623, 5879.
\apvalues AP:D5/D6  0,  0, 0, 0,  0, 0, 2,  11, 75, 533,  4493, 40906, 415069.
\end{countitem}

For orientable embeddings of dipoles, reversing the cyclic ordering at either vertex is equivalent to taking the \emph{Petrie dual} of the embedding, which twists all the edges (see \cite[Section 1.3]{E-MM13}).  We have two situations where our involution is Petrie duality.

\begin{countitem}(D3/D5)
For vertex-labeled orientable $n$-edge dipole embeddings
the number of Petrie-self-dual ones (S) is
$2\tx(D5)-\tx(D3) = \Adx(R_1)$,
and the number of non-Petrie-self-dual pairs (AP) is
$\tx(D3)-\tx(D5) = \rc4 (\Adx(I) + \Adx(R) - 2\Adx(R_1))$.

The number of Petrie-self-dual ones (S) for even positive $n$ appears in the OEIS \cite{OEIS} as A002866.

\svalues S:D3/D5   1,  1, 1, 0,  1, 0, 4,  0, 24, 0,  192, 0, 1920.
\apvalues AP:D3/D5  0,  0, 0, 1,  1, 4, 8,  39, 178, 1219,  9276, 83435, 834097.
\end{countitem}

\begin{countitem}(D4/D6)
For orientable $n$-edge dipole embeddings
the number of Petrie-self-dual ones (S) is
$2\tx(D6)-\tx(D4) = \rc2 (\Adx(R_1) + \Adx(R_1X))$,
and the number of non-Petrie-self-dual pairs (AP) is
$\tx(D6)-\tx(D4) = \rc8( \Adx(I) + \Adx(R) + 2\Adx(X)
			- 2\Adx(R_1) - 2\Adx(R_1X) )$.

\svalues S:D4/D6   1,  1, 1, 0,  1, 1, 3,  0, 15, 6,  102, 0, 990.
\apvalues AP:D4/D6  0,  0, 0, 1,  1, 3, 7,  28, 112, 680,  4873, 42529, 419958.
\end{countitem}

The next two situations involve formula (D7).
To discuss these, it is most natural to reinterpret (D1), (D5) and (D7) as counting hamilton cycles in a directed or undirected complete graph whose vertices are unlabeled but have a directed or undirected cyclic ordering.
In particular, (D1) counts directed hamilton cycles on a set of unlabeled vertices with a (directed) cyclic ordering, (D5) counts undirected hamilton cycles on a set of unlabeled vertices with a (directed) cyclic ordering and also directed hamilton cycles on a set of unlabeled vertices with an undirected cyclic ordering, and (D7) counts undirected hamilton cycles on a set of unlabeled vertices with an undirected cyclic ordering.

In this setting there are two natural involutions: reversal of a directed cyclic ordering, which may be considered a reflection, and arc-reversal of a directed cycle.  The quantities we obtain can be interpreted in terms of either of these involutions.

\begin{countitem}(D1/D7)
For directed hamilton cycles on a set of (unlabeled) vertices with a (directed) cyclic ordering
the number of reflexible ones (S) is
$2\tx(D7)-\tx(D1) = \Adx(R_1)$ (equal to (D3/D5:S), see values above),
and the number of chiral pairs (AP) is
$\tx(D1)-\tx(D7) = \rc2 ( \Adx(I) - \Adx(R_1) )$.

We may also interpret (S) as the number of arc-reversible ones and (AP) as the number of arc-irreversible pairs.

\smallskip
\apvalues AP:D1/D7  0,  0, 0, 1,  1, 4, 10,  54, 308, 2246,  18072, 164950, 1662434.
\end{countitem}

\begin{countitem}(D7/D5)
For undirected hamilton cycles on a set of (unlabeled) vertices with a (directed) cyclic ordering
the number of reflexible ones (S) is
$2\tx(D5)-\tx(D7) = \rc2 ( \Adx(R) + \Adx(R_1) )$,
and the number of chiral pairs (AP) is
$\tx(D7)-\tx(D5) = \rc4 ( \Adx(I) - \Adx(R) )$ (equal to (D2/D4:AP),
see values above).

We may also interpret these in terms of directed hamilton cycles on a set of (unlabeled) vertices with an undirected cyclic ordering: (S) is the number of arc-reversible ones, and (AP) is the number of arc-irreversible pairs.

\svalues S:D7/D5   1,  1, 1, 1,  2, 4, 10,  24, 72, 192,  672, 1920, 7680.
\end{countitem}

Our final two situations involve formula (D8).  To discuss these we re-interpret formulas (D3) and (D6) in terms of permutation matrices. 
Formula (D3) counts the equivalence classes of $\nn$ permutation matrices under cyclic shifts of the row set, cyclic shifts of the column set, and rotation by $180^\circ$.  
Formula (D6) counts the equivalence classes of $\nn$ permutation matrices under cyclic shifts of the row set, cyclic shifts of the column set, rotations of multiples of $90^\circ$, and transposition.

\begin{countitem}(D3/D8)
For equivalence classes of permutation matrices under cyclic shifts of row and column sets and rotations of $180^\circ$,
the number invariant under rotations of $90^\circ$ (S) is
$2\tx(D8)-\tx(D3) = \Adx(R_1 X)$,
and the number of pairs that are swapped by rotations of $90^\circ$ (AP) is
$\tx(D3)-\tx(D8) = \rc4 ( \Adx(I) + \Adx(R) - 2\Adx(R_1 X) )$.

\svalues S:D3/D8   1,  1, 1, 0,  1, 2, 2,  0, 6, 12,  12, 0, 60.
\apvalues AP:D3/D8  0,  0, 0, 1,  1, 3, 9,  39, 187, 1213,  9366, 83435, 835027.
\end{countitem}

\begin{countitem}(D8/D6)
For equivalence classes of permutation matrices under cyclic shifts of row and column sets and rotations of multiples of $90^\circ$,
the number that are symmetric (invariant under transposition) (S) is
$2\tx(D6)-\tx(D8) = \rc2 ( \Adx(R_1) + \Adx(X) )$,
and the number of asymmetric pairs (swapped by transposition) (AP) is (AP) is
$\tx(D8)-\tx(D6) = \rc8 ( \Adx(I) + \Adx(R) + 2\Adx(R_1 X)
			- 2\Adx(R_1) - 2\Adx(X) )$.

\svalues S:D8/D6   1,  1, 1, 1,  2, 3, 9,  17, 61, 147,  572, 1623, 6809.
\apvalues AP:D8/D6  0,  0, 0, 0,  0, 1, 1,  11, 66, 539,  4403, 40906, 414139.
\end{countitem}

\subsection{Bouquet results}\label{res-bouquet}

\subsubsection{Bouquet coset averages}\label{ca-res-bouquet}

There are two basic quantities for counting bouquet embeddings and related objects.
As we show in Part II, both have integer values, given positive integer inputs.
For $n=0$ and an arbitrary value of $k$ both of these quantities should be taken to be $1$; the following formulas apply for $n \ge 1$.
\begin{align*}
\Abx(I)
 &= \frac{1}{2n} \sum_{\substack{(d,g)\,:\, dg = 2n \\ \text{$d$ odd}}}
		\eul(d)\, \ma(g, g/2)\, d^{g/2}\, k^{g/2}
    + \frac{1}{2n} \sum_{\substack{(d,g)\,:\, dg = 2n \\ \text{$d$
even}}}
	 \eul(d) \sum_{j=0}^{\lfloor g/2 \rfloor} \ma(g, j)\, d^j\, k^{g-j};
\displaybreak[0]
\\[4pt]
 \Abx(R)
    &= \frac{1}{2} \left(
   \sum_{j=0}^{\lfloor n/2 \rfloor} \ma(n,j)\, 2^j\, k^{n-j}
   + \sum_{j=0}^{\lfloor (n-1)/2 \rfloor} \ma(n-1,j)\, 2^j\, k^{n-j}.
 \right).
\end{align*}
The values of $\Ab(I)(n,1)$ and $\Ab(R)(n,1)$ appear in the OEIS
\cite{OEIS}; see (B1) and (B1/B2:S) below.

\subsubsection{Counting basic bouquet objects}\label{bo-res-bouquet}

\begin{countitem}
(B1) The number of oriented embeddings of $k$-colored $n$-edge bouquets is $\Abx(I)$.

This is also the number of $k$-colored $n$-chord diagrams up to rotations (cyclic shifts).
This sequence for $k=1$ appears in the OEIS \cite{OEIS} as A007769.  It was also presented by Feng, Kwak, and Zhou \cite[Theorem 3.2]{FKZ10} specifically in the context of embeddings of bouquets.

\valuesk 1B1   1,  1, 2, 5,  18, 105, 902,  9749, 127072, 1915951,  32743182, 624999093, 13176573910.
\end{countitem}

\begin{countitem}
(B2) The number of orientable embeddings of $k$-colored $n$-edge bouquets is $\rc2 (\Abx(I)+\Abx(R))$.

This is also the number of $k$-colored $n$-chord diagrams up to rotations and reflections.
This sequence for $k=1$ appears in the OEIS \cite{OEIS} as A054499.

\valuesk 1B2   1,  1, 2, 5,  17, 79, 554,  5283, 65346, 966156,  16411700, 312700297, 6589356711.
\end{countitem}

\begin{countitem}
(B3) The number of generic (orientable or nonorientable) embeddings of $k$-colored $n$-edge bouquets is
$\rc2 (\Ab(I)(n, 2k)+\Ab(R)(n,2k))$
For $k=1$ this gives 
$$\rc2 (\Ab(I)(n, 2)+\Ab(R)(n,2))$$
as the number of generic embeddings of $n$-edge bouquets.

This sequence for $k=1$ was found by Kim and Park \cite[Theorem 3.2]{KP00}.  Their approach involves subdividing a bouquet $B_n$ to give a graph with $n$ triangles meeting at a common vertex.

\valuesk 1B3   1,  2, 6, 26,  173, 1844, 29570,  628680, 16286084, 490560202,  16764409276, 639992710196, 26985505589784.
\end{countitem}

\begin{countitem}
(B4) The number of nonorientable embeddings of $k$-colored $n$-edge bouquets is
$\text{(B3)}-\text{(B2)} =
	\rc2 (\Ab(I)(n,2k)+\Ab(R)(n,2k)-\Abx(I)-\Abx(R))$.
For $k=1$ this gives
	$$\rc2 (\Ab(I)(n,2)+\Ab(R)(n,2)-\Ab(I)(n,1)-\Ab(R)(n,1))$$
as the number of nonorientable embeddings of $n$-edge bouquets.

\valuesk 1B4   0,  1, 4, 21,  156, 1765, 29016,  623397, 16220738, 489594046,  16747997576, 639680009899, 26978916233073.
\end{countitem}

\subsubsection{Symmetric and asymmetric bouquet objects}\label{sa-res-bouquet}

Again we provide some applications of equations \eqref{symmetric} and \eqref{asymmetric} of Subsection \ref{term-inv}.  Given an involution (here, just reflection) that acts on objects counted by (B$i$) and creates equivalence classes counted by (B$j$), item  (B$i$/B$j$) states the formulas for the numbers of symmetric objects, (S) or (B$i$/B$j$:S), and asymmetric pairs of objects, (AP) or (B$i$/B$j$:AP).

\begin{countitem}
(B1/B2) For $n$-edge bouquet embeddings
the number of reflexible ones (S) is
$2\text{(B2)} - \text{(B1)} = \Abx(R)$, and
the number of chiral pairs (AP) is
$\text{(B1)} - \text{(B2)} = \rc2 (\Abx(I)-\Abx(R))$.

The number of reflexible ones (S) for $k=1$ appears in the OEIS \cite{OEIS} as A018191.  It was also presented by Feng, Kwak, and Zhou \cite[Theorem 4.2]{FKZ13} specifically in the context of embeddings of bouquets.
The number of chiral pairs (AP) for $k=1$ occurs in the OEIS \cite{OEIS} as A054938.

\svaluesk 1S:B1/B2   1,  1, 2, 5,  16, 53, 206,  817, 3620, 16361,  80218, 401501, 2139512.
\apvaluesk 1AP:B1/B2  0,  0, 0, 0,  1, 26, 348,  4466, 61726, 949795,  16331482, 312298796, 6587217199.
\end{countitem}

\subsection{Directed bouquet results}\label{res-dirbouquet}

\subsubsection{Directed bouquet coset averages}\label{ca-res-dirbouquet}

There are four basic quantities for counting directed embeddings of directed bouquets and related objects.
As we show in Part II, all have integer values, given positive integer inputs.
For $n=0$ and an arbitrary value of $k$ all of these quantities should be taken to be $1$; the following formulas apply for $n \ge 1$.
\begin{align*}
 \Aax(I)
	&= \frac{1}{n} \sum_{(d,g): dg = n} \phi(d)\, g!\, d^g\, k^g ;
\displaybreak[0]
\\[4pt]
 \Aax(R) &= \begin{cases}
 0 & \text{if $n$ is even,} \\
 \displaystyle \left(\frac{n-1}{2}\right)!\, 2^{(n-1)/2}\, k^{(n+1)/2}
   & \text{if $n$ is odd;} \\
\end{cases}
\displaybreak[0]
\\[4pt]
 \Aax(F) 
  &= \begin{cases}
	0 & \text{if $n$ is even,} \\
   \displaystyle \frac{1}{n}
	\sum_{\substack{(d,g)\,:\, dg = n \\ \text{$g$ odd}}}
	\eul(2d) \sum_{j=0}^{\lfloor g/2 \rfloor}
	\ma(g, j)\, d^j\, k^{g-j}
	& \text{if $n$ is odd;} \\
   \end{cases}
\displaybreak[0]
\\[4pt]
 \Aax(RF)
     &= \sum_{j=1}^{\lfloor n/2 \rfloor} \ma(n,j)\, k^{n-j} .
\end{align*}
The values of $\Aa(I)(n,1)$, $\Aa(R)(n,1)$, and $\Aa(RF)(n,1)$ appear in the OEIS; see (A1), (A1/A2:S), and (A1/A5:S) below.

\subsubsection{Counting basic directed bouquet objects}\label{bo-res-dirbouquet}

\begin{countitem}
(A1) The number of oriented $k$-colored directed embeddings of $n$-arc directed bouquets is $\Aax(I)$.

This sequence for $k=1$ appears in the OEIS \cite{OEIS} as A061417.
It was also presented by Chen, Gao, and Huang \cite[Theorem 3.3]{CGH18} specifically in the context of directed embeddings of directed bouquets.

\valuesk 1A1   1,  1, 2, 4,  10, 28, 136,  726, 5100, 40362,  363288, 3628810, 39921044.
\end{countitem}

\begin{countitem}
(A2) The number of orientable $k$-colored directed embeddings of $n$-arc directed bouquets is
$\rc2 (\Aax(I)+\Aax(R))$.

\valuesk 1A2   1,  1, 1, 3,  5, 18, 68,  387, 2550, 20373,  181644, 1816325, 19960522.
\end{countitem}

\begin{countitem}
(A3) The number of arc-reversal classes of oriented $k$-colored directed embeddings of $n$-arc directed bouquets is
$\rc2 (\Aax(I)+\Aax(F))$.

\valuesk 1A3   1,  1, 1, 3,  5, 17, 68,  380, 2550, 20328,  181644, 1816028, 19960522.
\end{countitem}

\begin{countitem}
(A4) The number of arc-reversal classes of orientable $k$-colored directed embeddings of $n$-arc directed bouquets is
$\rc4 (\Aax(I) + \Aax(R) + \Aax(F) + \Aax(RF))$.

\valuesk 1A4   1,  1, 1, 3,  5, 17, 53,  260, 1466, 10915,  93196, 917898, 10015299.
\end{countitem}

\begin{countitem}
(A5) The number of classes of oriented $k$-colored directed embeddings of $n$-arc directed bouquets under simultaneous reflection and arc reversal is
$\rc2 ( \Aax(I) + \Aax(RF) )$.

\valuesk 1A5   1,  1, 2, 4,  10, 27, 106,  479, 2932, 21491,  186392, 1832253, 20030598 .
\end{countitem}

\begin{countitem}
(A6) The number of generic (orientable or nonorientable) $k$-colored directed embeddings of $n$-arc directed bouquets is
$\rc2 ( \Aa(I)(n,2k)+\Aa(R)(n,2k) )$.
For $k=1$ this gives
$$\rc2 ( \Aa(I)(n,2)+\Aa(R)(n,2) )$$
as the number of generic directed embeddings of directed bouquets.

\valuesk 1A6   1,  2, 3, 14,  54, 420, 3886,  46470, 645524, 10328214, 
185800748, 3716014090, 81749732156.
\end{countitem}

\begin{countitem}
(A7) The number of nonorientable $k$-colored directed embeddings of $n$-arc directed bouquets is
$\text{(A6)} - \text{(A2)} =
	\rc2 ( \Aa(I)(n,2k)+\Aa(R)(n,2k) - \Aa(I)(n,k)-\Aa(R)(n,k) )$.
For $k=1$, this gives
$$\rc2 ( \Aa(I)(n,2)+\Aa(R)(n,2) - \Aa(I)(n,1)-\Aa(R)(n,1) ) $$
as the number of nonorientable directed embeddings of directed bouquets.

\valuesk 1A7   0,  1, 2, 11,  49, 402, 3818,  46083, 642974, 10307841,  185619104, 3714197765, 81729771634.
\end{countitem}

\begin{countitem}
(A8) The number of arc-reversal classes of generic (orientable or nonorientable) $k$-colored directed embeddings of $n$-arc directed bouquets is
$\rc4 (\Aa(I)(n,2k) + \Aa(R)(n,2k) + \Aa(F)(n,2k) + \Aa(RF)(n,2k))$.
For $k=1$ this gives
$$\rc4 (\Aa(I)(n,2) + \Aa(R)(n,2) + \Aa(F)(n,2) + \Aa(RF)(n,2))$$
as the number of arc-reversal classes of generic directed embeddings of directed bouquets.

\valuesk 1A8   1,  2, 3, 14,  46, 304, 2289,  25096, 330862, 5211052,  93130670, 1859431284, 40882543694.
\end{countitem}

\begin{countitem}
(A9) The number of arc-reversal classes of nonorientable $k$-colored directed embeddings of $n$-arc directed bouquets is
$\text{(A8)} - \text{(A4)} =
 \rc4 (\Aa(I)(n,2k) + \Aa(R)(n,2k) + \Aa(F)(n,2k) + \Aa(RF)(n,2k)
  -\Aa(I)(n,k) - \Aa(R)(n,k) - \Aa(F)(n,k) - \Aa(RF)(n,k))$.
For $k=1$ this gives
$$ \rc4 (\Aa(I)(n,2) + \Aa(R)(n,2) + \Aa(F)(n,2) + \Aa(RF)(n,2)
  -\Aa(I)(n,1) - \Aa(R)(n,1) - \Aa(F)(n,1) - \Aa(RF)(n,1))$$
as the number of arc-reversal classes of nonorientable directed embeddings of directed bouquets.

\valuesk 1A9   0,  1, 2, 11,  41, 287, 2236,  24836, 329396, 5200137,  93037474, 1858513386, 40872528395.
\end{countitem}

\subsubsection{Symmetric and asymmetric directed bouquet objects}\label{sa-res-dirbouquet}

Again we provide some applications of equations \eqref{symmetric} and \eqref{asymmetric} of Subsection \ref{term-inv}. 
Given an involution that acts on objects counted by (A$i$) and creates equivalence classes counted by (A$j$), item (A$i$/A$j$) states the formulas for the numbers of symmetric objects, (S) or (A$i$/A$j$:S), and asymmetric pairs of objects, (AP) or (A$i$/A$j$:AP).

First we consider situations where our involution is `reflection' in the sense of reversal of the surface orientation.

\begin{countitem}
(A1/A2) For oriented $k$-colored directed embeddings of $n$-arc directed bouquets
the number of reflexible ones (S) is
$2\text{(A2)}-\text{(A1)} = \Aax(R)$, and
the number of chiral pairs (AP) is
$\text{(A1)}-\text{(A2)} = \rc2 (\Aax(I)-\Aax(R))$.

The number of reflexible ones for $k=1$ and odd $n$ appears in the OEIS \cite{OEIS} as A000165.

\svaluesk  1S:A1/A2   1,  1, 0, 2,  0, 8, 0,  48, 0, 384,  0, 3840, 0.
\apvaluesk  1AP:A1/A2  0,  0, 1, 1,  5, 10, 68,  339, 2550, 19989,  181644, 1812485, 19960522.
\end{countitem}

\begin{countitem}
(A3/A4) For arc-reversal classes of $k$-colored directed embeddings of $n$-arc directed bouquets
the number of reflexible ones (S) is
$2\text{(A4)}-\text{(A3)} = \rc2 (\Aax(R) + \Aax(RF))$, and
the number of chiral pairs (AP) is
$\text{(A3)}-\text{(A4)}
	= \rc4 (\Aax(I)+\Aax(F)-\Aax(R)-\Aax(RF))$.

\svaluesk  1S:A3/A4   1,  1, 1, 3,  5, 17, 38,  140, 382, 1502,  4748, 19768, 70076.
\apvaluesk  1AP:A3/A4  0,  0, 0, 0,  0, 0, 15,  120, 1084, 9413,  88448, 898130, 9945223.
\end{countitem}

Next we consider situations where our involution is arc-reversal.

\begin{countitem}
(A1/A3) For oriented $k$-colored directed embeddings of $n$-arc directed bouquets
the number of arc-reversible ones (S) is
$2\text{(A3)}-\text{(A1)} = \Aax(F)$, and
the number of arc-irreversible pairs (AP) is
$\text{(A1)}-\text{(A3)} = \rc2 (\Aax(I)-\Aax(F))$.

\svaluesk  1S:A1/A3   1,  1, 0, 2,  0, 6, 0,  34, 0, 294,  0, 3246, 0.
\apvaluesk  1AP:A1/A3  0,  0, 1, 1,  5, 11, 68,  346, 2550, 20034,  181644, 1812782, 19960522.
\end{countitem}

\begin{countitem}
(A2/A4) For orientable $k$-colored directed embeddings of $n$-arc directed bouquets
the number of arc-reversible ones (S) is
$2\text{(A4)}-\text{(A2)} = \rc2 (\Aax(F)+\Aax(RF))$, and
the number of arc-irreversible pairs (AP) is
$\text{(A2)}-\text{(A4)} =
	\rc4 (\Aax(I)+\Aax(R)-\Aax(F)-\Aax(RF))$.

\svaluesk  1S:A2/A4   1,  1, 1, 3,  5, 16, 38,  133, 382, 1457,  4748, 19471, 70076.
\apvaluesk  1AP:A2/A4  0,  0, 0, 0,  0, 1, 15,  127, 1084, 9458,  88448, 898427, 9945223.
\end{countitem}

\begin{countitem}
(A6/A8) For generic $k$-colored directed embeddings of $n$-arc directed bouquets
the number of arc-reversible ones (S) is
$2\text{(A8)}-\text{(A6)} =
	\rc2 (\Aa(F)(n,2k)+\Aa(RF)(n,2k))$, and
the number of arc-irreversible pairs (AP) is
$\text{(A6)}-\text{(A8)} =
	\rc4 (\Aa(I)(n,2k)+\Aa(R)(n,2k)-\Aa(F)(n,2k)-\Aa(RF)(n,2k))$.

\svaluesk  1S:A6/A8   1,  2, 3, 14,  38, 188, 692,  3722, 16200, 93890,  460592, 2848478, 15355232.
\apvaluesk  1AP:A6/A8  0,  0, 0, 0,  8, 116, 1597,  21374, 314662, 5117162,  92670078, 1856582806, 40867188462.
\end{countitem}

\begin{countitem}
(A7/A9) For nonorientable $k$-colored directed embeddings of $n$-arc directed bouquets
the number of arc-reversible ones (S) is
$2\text{(A9)}-\text{(A7)} =
	\rc2 (\Aa(F)(n,2k)+\Aa(RF)(n,2k)-\Aa(F)(n,k)-\Aa(RF)(n,k))$, and
the number of arc-irreversible pairs (AP) is
$\text{(A7)}-\text{(A9)} =
	\rc4 (\Aa(I)(n,2k)+\Aa(R)(n,2k)-\Aa(F)(n,2k)-\Aa(RF)(n,2k)
	-\Aa(I)(n,k)-\Aa(R)(n,k)+\Aa(F)(n,k)+\Aa(RF)(n,k))$.

\svaluesk  1S:A7/A9   0,  1, 2, 11,  33, 172, 654,  3589, 15818, 92433,  455844, 2829007, 15285156.
\apvaluesk  1AP:A7/A9  0,  0, 0, 0,  8, 115, 1582,  21247, 313578, 5107704,  92581630, 1855684379, 40857243239.
\end{countitem}

Finally we consider situations involving simultaneous reflection and arc reversal.

\begin{countitem}
(A1/A5) For oriented $k$-colored directed embeddings of $n$-arc directed bouquets
the number symmetric under simultaneous reflection and arc reversal (S) is
$2\tx(A5)-\tx(A1) = \Aax(RF)$,
and the number of asymmetric pairs (AP) is
$\tx(A1)-\tx(A5) = \rc2 ( \Aax(I) - \Aax(RF) )$.

The number of symmetric ones (S) for $k=1$ appears in the OEIS \cite{OEIS} as A000085; it is the total number of matchings in an $n$-vertex complete graph.

\svaluesk 1S:A1/A5   1,  1, 2, 4,  10, 26, 76,  232, 764, 2620,  9496, 35696, 140152.
\apvaluesk 1AP:A1/A5  0,  0, 0, 0,  0, 1, 30,  247, 2168, 18871, 176896, 1796557, 19890446.
\end{countitem}

\begin{countitem}
(A5/A4) Let $\cS$ be the set of equivalence classes of oriented $k$-colored directed embeddings of $n$-arc directed bouquets under simultaneous reflection and arc reversal. Considering elements of $\cS$ up to reflection is the same as considering elements of $\cS$ up to arc-reversal.
The number of reflexible (or arc-reversible) elements of $\cS$ (S) is
$2\tx(A4)-\tx(A5) = \rc2 ( \Aax(R) + \Aax(F) )$,
and the number of chiral (or arc-irreversible) pairs (AP) is
$\tx(A5)-\tx(A4) = \rc4 ( \Aax(I) + \Aax(RF) - \Aax(R) - \Aax(F) )$.

\svaluesk 1S:A5/A4   1,  1, 0, 2,  0, 7, 0,  41, 0, 339,  0, 3543, 0.
\apvaluesk 1AP:A5/A4  0,  0, 1, 1,  5, 10, 53,  219, 1466, 10576,  93196, 914355, 10015299.
\end{countitem}

\section{Future directions}

There are a number of problems that follow naturally from the work in this paper.  For bouquets and directed bouquets we have found the number of generic (orientable or nonorientable) and nonorientable embeddings.  It would also be natural to count equivalence classes of generic embeddings under Petrie duality.  For dipoles we have not yet counted generic and nonorientable embeddings, and counting those would be a natural next step, building on some of the ideas we used in counting dipolar cogs.  We have already counted Petrie duality classes of orientable dipole embeddings, but it would also be interesting  do this for generic dipole embeddings. Chen, Gao, and Huang \cite[Theorem 2.4]{CGH18} have counted directed embeddings of Eulerian directed dipoles in oriented surfaces, and it should be possible to obtain related results similar to the results in Subsections \ref{res-dipole} and \ref{res-dirbouquet}. All of these problems seem approachable by extending the techniques used here.

Another open problem is to enumerate the looped dipoles of Figure \ref{fig:Upper2}, which is equivalent  to counting upper embeddable graphs with two faces.  A similar problem is to count 'pointed' graphs, that is, graphs with one distinguished vertex that is incident to every edge.  This vertex is then incident to a collection of loops and digons as in Figure \ref{fig:EdgeOuter}.  Such graphs are dual to edge-outer embeddable graphs.
Because these problems involve multiple graphs for a given number of edges, they seem more difficult than the problems discussed in the first paragraph. However, the formulas and approaches given here can likely serve as a foundation for further work on looped dipoles and pointed graphs.

\part{\Large Technical details}
This second part of the paper provides the technical details and formal proofs of the results in the first part of the paper.
\section{Basic counting results}

We begin with some basic counting results we will use. The first is very well known.

\begin{thm}[Burnside's Lemma (stated earlier by Cauchy and Frobenius)]\label{burnside}
Suppose $\Ga$ is a group acting on a set $S$.  Then the number of
orbits of the action, i.e., the number of equivalence classes under the
symmetries provided by $\Ga$, is
$$ \frac{1}{|\Ga|} \sum_{\ga \in \Ga} |\fix(\ga)|, $$
where $\fix(\ga)$ is the set of elements of $S$ fixed by $\ga$.
\end{thm}

In several places we need to count how many permutations $\tau$ have $\tau^2 = \al$ for a given permutation $\al$.  We will use the following lemma.  Recall that $\ma(n,j)$ is the number of $j$-matchings of an $n$-set.


\begin{lemma}\label{squareroots}
For integers $\ell \ge 1$ and $m \ge 0$ define
$$\sr(\ell, m) = \begin{cases}
   0	& \text{if $\ell$ is even and $m$ is odd,} \\[2pt]
   \ma(m, m/2)\, \ell^{m/2}
	& \text{if $\ell$ is even and $m$ is even,} \\[2pt]
   \displaystyle\sum_{j=0}^{\lfloor m/2 \rfloor} \ma(m, j)\, \ell^j
	& \text{if $\ell$ is odd.}
 \end{cases}
$$
Notice that $\sr(\ell,m) = 1$ if $m=0$.
Let $\al$ be a permutation of an $n$-set with $a_\ell$ cycles of length $\ell$ for each $\ell$, $1 \le \ell \le n$.  Then the number of permutations $\tau$ with $\tau^2 = \al$ is $\displaystyle\prod_{\ell=1}^n \sr(\ell, a_\ell)$.
\end{lemma}

\begin{proof}
We need to consider how the cycles of $\al$ come from the cycles of $\tau$.
A cycle of even length $2\ell$ in $\tau$ yields two cycles of length $\ell$ in $\al$.  A cycle of odd length in $\tau$ yields a cycle of the same odd length in $\al$.  Therefore, when $\ell$ is even each cycle of length $\ell$ in $\al$ must be paired with another such cycle and come from a cycle of length $2\ell$ in $\tau$.  When $\ell$ is odd each cycle of length $\ell$ in $\al$ can either come from a cycle of length $\ell$ in $\tau$, or be paired with another cycle and come from a cycle of length $2\ell$ in $\tau$.

If $\al$ has an odd number $m=a_\ell$ of cycles of even length $\ell$, we cannot pair them all up, so the number of possibilities for the corresponding cycles in $\tau$ is $0$, which is $\sr(\ell,m)=\sr(\ell, a_\ell)$.

If $\al$ has an even number $m=a_\ell$ of cycles of even length $\ell$, then there are $\ma(m, m/2)$ ways to pair these up, and for each pairing there are $\ell$ ways to interleave two cycles of length $\ell$ to obtain a cycle of length $2\ell$.  Therefore, the number of possibilities in $\tau$ is $\ma(m, m/2) \ell^{m/2}$, which is $\sr(\ell,m) = \sr(\ell, a_\ell)$.

If $\al$ has $m=a_\ell$ cycles of odd length $\ell$, then for each $j$ with $0 \le j \le \lfloor m/2 \rfloor$ we can pair up $2j$ of these cycles in $\ma(m, j)$ ways, interleave each of the $j$ pairs in $\ell$ ways, and use the unique square root of each of the $m-2j$ unpaired cycles, giving giving $\ma(m, j) \ell^j$ possibilities.
Summing over all $j$ gives that the number of possibilities for the corresponding cycles in $\tau$ is
$\sum_{j=0}^{\lfloor m/2 \rfloor} \ma(m, j) \ell^j$,
which is $\sr(\ell,m) = \sr(\ell, a_\ell)$.

Now multiplying the number of possibilities for each $\ell$ gives the result.
\end{proof}

We use the following easy argument in several places so it is convenient to summarize it here for reference.

\begin{observation}\label{matchingpairs}
Suppose we have sets $A, B$ with $|A|=|B|=2k$, and partitions $\mathcal{A}, \mathcal{B}$ of $A$ and $B$, respectively, into pairs ($2$-subsets).  Let $\psi : A \to B$ be a bijection that preserves pairs, i.e., such that $\psi(A') \in \mathcal{B}$ for all $A' \in \mathcal{A}$.
There are $k!$ choices for the bijection that $\psi$ establishes a between $\mathcal{A}$ and $\mathcal{B}$, and $2$ ways for $\psi$ to map the elements of each $A'$ to the elements of $\psi(A')$.  Hence the number of possible maps $\psi$ is $k!\, 2^k$.
\end{observation}

\section{Proofs for dipole formulas}\label{pf-dipole}

\subsection{Labeled dipoles and symmetry operations}\label{so-pf-dipole}

In this section we prove the counting results from Subsection \ref{res-dipole} regarding embeddings of dipoles and related objects.  Recall that a dipole $D_n$ has two vertices and $n$ edges, each edge having both vertices as its ends (so there are multiple edges but no loops).  We think of each edge as consisting of two half-edges, each incident with one of the vertices.

Our results on dipoles will be proved by elementary techniques (straightforward applications of Burnside's Lemma) based on groups acting on a set of objects that we will call labeled dipoles. As we will see, labeled dipoles are in one-to-one correspondence with elements of the symmetric group $\Sz_n$, so our results can also be interpreted as counting results for permutations, or permutation matrices, under various equivalence relations.

A \emph{labeled dipole} is a dipole $D$ where the vertices receive distinct labels $0$ and $1$, and for each vertex the half-edges incident with that vertex receive distinct labels from $\mZ_n = \{0, 1, 2, \dots, n-1\}$, where $n = |E(D)|$.
An example was given in Figure \ref{fig:dipolesII}.
We let $\cD_n$ denote the set of $n$-edge labeled dipoles. For $j \in \{0,1\}$, we refer to the vertex labeled $j$ as \emph{vertex $j$}, and the half-edges incident to vertex $j$ as \emph{$j$-half-edges}.
For each edge $e$ we let $\lab_j(e)$ denote the label of the $j$-half-edge of $e$.

A labeled dipole $D$ is completely described by the set of ordered pairs $\{(\lab_0(e), \lab_1(e)) \;|\allowbreak\; e \in E(D)\}$.  There is therefore a one-to-one correspondence between $n$-edge labeled dipoles and sets
$P = \{(a_0, b_0), (a_1, b_1), \dots, (a_{n-1}, b_{n-1})\}$ with
$\{a_0, a_1, \dots, a_{n-1}\} = \{b_0, b_1, \dots, b_{n-1}\} \allowbreak = \mZ_n$.
Let $\cP_n$ denote the collection of all such sets $P$.
Recall that a function is formally defined as a set of ordered pairs whose first components are distinct.
Therefore, $P$ may be regarded as the formal representation of a function $\pi: \mZ_n \to \mZ_n$ with $\pi(a_i) = b_i$ for $0 \le i \le n-1$.  It is easy to see that $\pi$ is a bijection. Thus, $\pi$ is a permutation of $\mZ_n$, an element of the symmetric group $\symz_n$.  Thus, we have natural bijections between $\cD_n$, $\cP_n$ and $\symz_n$.

For example, the labeled dipole in Figure \ref{fig:dipolesII} corresponds to $\{(0,4), (1,2), (2,3),\allowbreak (3,0),\allowbreak (4,1)\} \in \cP_5$ and the permutation $
\left[ \begin{matrix}
  0 & 1 & 2 & 3 & 4 \\ 4 & 2 & 3 & 0 & 1 \\
\end{matrix} \right] \in \symz_5$.

Each object of types (D1)--(D8) above can be turned into a labeled dipole (or permutation of $\mZ_n$) in a natural way.  But some choices are involved in doing this, so each object can be turned into several different labeled dipoles (or permutations), which we wish to characterize as being equivalent under certain symmetries.  The symmetries always include cyclically shifting the labels of the $0$-half-edges and cyclically shifting the labels of the $1$-half-edges, and may also include reversing the labels of the $0$-half-edges or the $1$-half-edges, or swapping the vertex labels.

Our symmetries can be defined in terms of how they act on ordered pairs $(a,b) \in \pz_n$.  However, to simplify some proofs, and to avoid treating $n=1$ and $2$ as special cases, we will define them more generally, as permutations of the set $\pr_n$, where $\mR_n = \mR/n\mR$ is the additive group of real numbers modulo $n$, whose underlying set can be identified with the real interval $[0,n)$.
Note that $\mR_n$ contains $\mZ_n$ as a subgroup.
We define $S_0, S_1, R_0, R_1, X \in \sym(\mR_n \times \mR_n)$ as follows:

 \listitem $S_0(a, b) = (a+1, b)$;\quad
      $S_1(a, b) = (a, b+1)$;
 \listitem $R_0(a, b) = (-a, b)$;\quad
      $R_1(a, b) = (a, -b)$; \quad and
 \listitem $X(a, b) = (b, a)$.

\noindent
We let $\Wgd = \la S_0, S_1, R_0, R_1, X \ra \le \sym(\pr_n)$.

Consider the induced action of $\sym(\pr_n)$ on $P \subseteq \pr_n$ by $T(P) = \{T(a,b) \st (a,b) \in P\}$. If $T \in \{S_0, S_1, R_0, R_1, X\}$ then both $T$ and $T\iv$ map $\pz_n$ to itself, and moreover both $T$ and $T\iv$ map $\cP_n$ to itself under the induced subset action.
It follows that these statements also hold for all $T \in \Wgd$.  So we have an action of $\Wgd$ on $\cP_n$, and using the natural bijections between $\cP_n$ and $\cD_n$ or $\symz_n$, we also obtain actions of $\Wgd$ on those sets.

The actions of $S_0$ and $S_1$ on $\cD_n$ correspond to shifting the labels of the $0$- and $1$-half-edges, respectively; $R_0$ and $R_1$ correspond to reversing (specifically, negating) the labels of the $0$- and $1$-half-edges, respectively, and $X$ corresponds to swapping the vertex labels.

We note that `reflections' of labels can be regarded as elements of a dihedral group generated by cyclic shifts (rotations) and one reflection.  In particular, for labels around vertex $i$ (assumed to be equally spaced around a circle), $S_i^h R$ maps $j \mapsto h-j$, which is reflection about an axis passing through $h/2$ and $h/2+n/2$.

In some contexts `reversal' of labels might naturally be regarded as the map on $\mZ_n$ with $0 \mapsto n-1$, $1 \mapsto n-2$, $\dots$, $n-1 \mapsto 0$.  This maps $j \mapsto n-1-j$, and is just reflection about an axis through $n-1/2$ or $(n-1)/2$, and for labels around vertex $i$ is represented by the transformation $V_i = S_i^{n-1} R_i$. However, using $R_0$ and $R_1$ instead of $V_0$ and $V_1$ generally simplifies calculations.
Since we will always include $S_0$ and $S_1$ in our groups, and since $\la S_j, R_j \ra = \la S_j, V_j \ra$ for $j \in \{0,1\}$, working with $R_0$ and $R_1$ is equivalent to working with $V_0$ and $V_1$, respectively.

The actions on $\Sz_n$ are as follows. Let $\si, \rho \in \Sz_n$ be defined by $\si(i) = i+1$ and $\rho(i) = -i$, and let $\pi$ be an arbitrary element of $\Sz_n$.
First, $S_0$ represents a cyclic shift in the input variable of $\pi$: instead of $i \mapsto \pi(i)$, we have $i+1 \mapsto \pi(i)$, or $j \mapsto \pi(j-1) = \pi\si\iv(j)$.
Next, $S_1$ represents a cyclic shift in the output variable: instead of $i \mapsto \pi(i)$, we have $i \mapsto \pi(i)+1 = \si \pi (i)$.
Now, $R_0$ represents reversal (negation) of the input variable: instead of $i \mapsto \pi(i)$, we have $-i \mapsto \pi(i)$, or $j \mapsto \pi(-j)=\pi\rho(j)$.
Next, $R_1$ represents reversal of the output variable: instead of $i \mapsto \pi(i)$, we have $i \mapsto -\pi(i) = \rho\pi(i)$.
Finally, $X$ corresponds to swapping the input and output variables, which inverts the permutation: instead of $i \mapsto \pi(i)$, we have $\pi(i) \mapsto i$, or $j \mapsto \pi\iv(j)$.  To summarize,
\begin{tightequation}
 S_0(\pi) = \pi \si\iv,\quad
 S_1(\pi) = \si \pi,\quad
 R_0(\pi) = \pi \rho,\quad
 R_1(\pi) = \rho\, \pi,\quad
 \text{ and }\quad
 X(\pi) = \pi\iv.
\end{tightequation}

%
\begin{figure}
\centering
\begin{subfigure}{0.7\textwidth}
    \includegraphics[clip, trim=2cm 20cm 0cm 2cm, width=\textwidth]{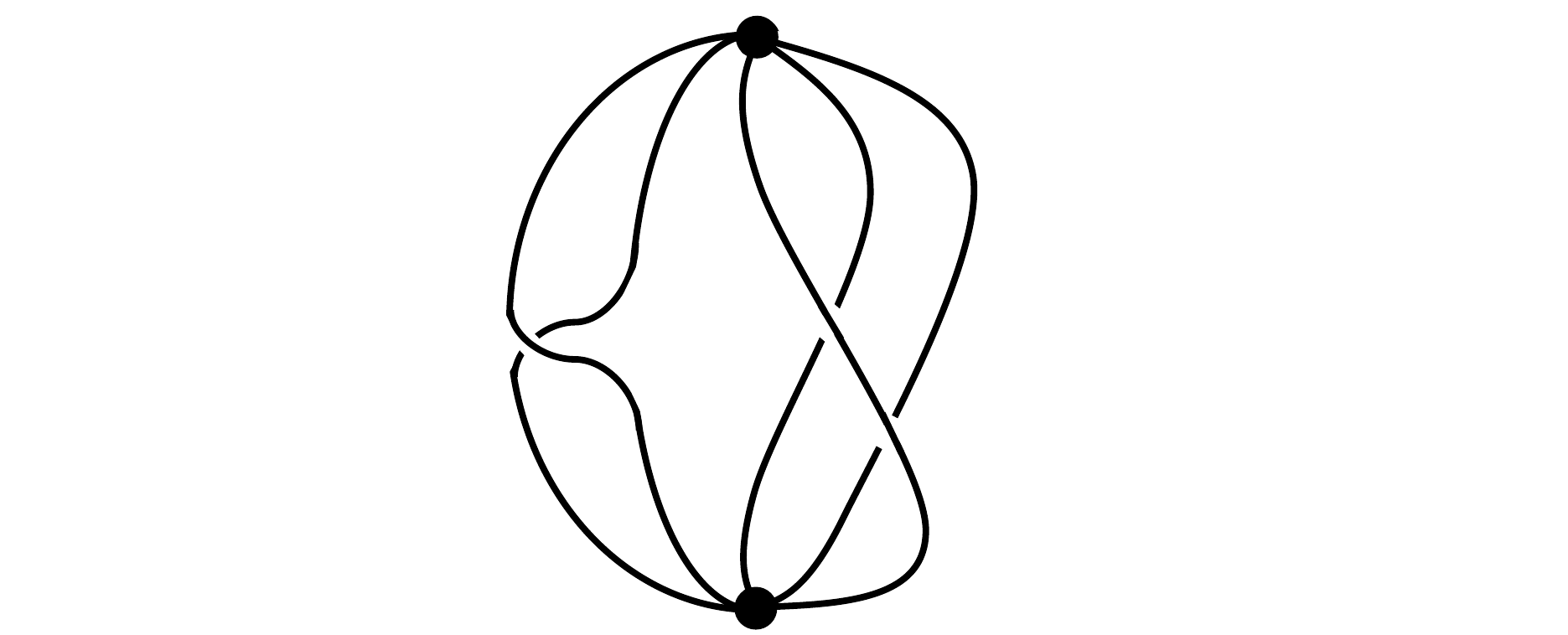}
    \caption{An unlabeled dipole embedded in an oriented surface.}
    \label{fig:dipoleblack}
\end{subfigure}
\begin{subfigure}{0.3\textwidth}
\centering
    \includegraphics[clip, trim=1cm 10cm 12cm 10cm, width=\textwidth]{figures/dipoleperms.pdf}
    \vskip10pt
    \caption{L1: An arbitrary vertex and edge labeling.}
    \label{fig:dipoleV1}
\end{subfigure}
\hspace{2cm}
\begin{subfigure}{0.3\textwidth}
    \includegraphics[clip, trim=12cm 9cm 1cm 10cm, width=\textwidth]{figures/dipoleperms.pdf}
    \caption{L2: A different arbitrary vertex and edge labeling.}
    \label{fig:dipoleV2}
\end{subfigure} 
\break\vskip20pt
\begin{subfigure}{0.3\textwidth}
    \includegraphics[clip, trim=1cm 0cm 11cm 20cm, width=\textwidth]{figures/dipoleperms.pdf}
    \caption{Extracting the permutation for L1.}
    \label{fig:permdipoleV1}
\end{subfigure}
\hspace{2cm}
\begin{subfigure}{0.3\textwidth}
    \includegraphics[clip, trim=9.8cm 0cm 2.2cm 20.4cm, width=\textwidth]{figures/dipoleperms.pdf}
    \caption{Extracting the permutation for L2.}
    \label{fig:permdipoleV2}
\end{subfigure}
\vskip3pt
\caption{There are many ways to label the vertices and edges of an embedded dipole, but any labelings of the same embedded dipole will correspond to permutations related by the symmetry operations.  }
\label{fig:dipoleperms}
\end{figure}

Figure \ref{fig:dipoleperms} illustrates how we may obtain two equivalent labeled dipoles from the oriented embedding of a vertex-unlabelled dipole in Figure \ref{fig:dipoleblack}, and describe the relationship using the symmetries we have just defined.
The permutation
$\pi_1= \left[\begin{matrix} 0&1&2&3&4\\3&2&4&0&1     \end{matrix}\right]$
corresponds to the labeling L1 in Figure \ref{fig:dipoleV1} and may be read from the diagram in Figure \ref{fig:permdipoleV1}.   Similarly, the permutation $\pi_2=\left[\begin{matrix} 0&1&2&3&4\\3&0&1&2&4     \end{matrix}\right]$ corresponds to the labeling L2 in Figure \ref{fig:dipoleV2} and may be read from the diagram in Figure \ref{fig:permdipoleV2}. The permutation $\pi_1$ is related to the permutation $\pi_2$ by exchanging the vertex labels (rotating the diagram by $180^{\circ}$), then incrementing the edge labels on the top by $2$ and those on the bottom by $3$.  Formally, this is:  
\begin{equation*}
\pi_2 =S^2_0 S^3_1 X (\pi_1) \quad\text{ or }\quad
	\pi_2 = \sigma^3 \pi_1^{-1} \sigma^{-2}.
\end{equation*}

The following theorem provides some relationships between the generators of $\Wgd$.

\begin{thm}\label{dsymmetries}
Let $n$ be a positive integer, and consider $S_0, S_1, R_0, R_1, X \in \sym(\mR_n \times \mR_n)$ as defined above.  Let $I$ be the identity of $\sym(\mR_n \times \mR_n)$.

\smallskip\noindent
(a) Then (composing functions right to left):

 \listitem $S_0^n = S_1^n = R_0^2 = R_1^2 = X^2 = I$;
 \listitem each of $S_0$ or $R_0$ commutes with each of $S_1$ or $R_1$;
 \listitem $R_0 S_0 = S_0\iv R_0$,\quad
      and\quad $R_1 S_1 = S_1\iv R_1$;
 \listitem $X S_0 = S_1 X$,\quad
      and\quad $X R_0 = R_1 X$.

\noindent
Moreover, if $R = R_0 R_1$ so that $R(a,b) = (-a, -b)$ then

 \listitem $R^2=I$,\quad
     $R S_0 = S_0\iv R$,\quad
     $R S_1 = S_1\iv R$,\quad
     and\quad $XR = RX$.

\smallskip\noindent
(b) Every element of $\Wgd = \la S_0, S_1, R_0, R_1, X \ra$ can be written uniquely as $S_0^h\, S_1^k\, R_0^p\, R_1^q\, X^s$ for some $h, k \in \mZ_n$ and $p, q, s \in \{0,1\}$.  Hence $|\Wgd| = 8n^2$.
\end{thm}

\begin{proof}
(a) All of these relations can be checked easily from the definitions.

(b) To get that each $T \in \Ga$ can be expressed as  $S_0^h\, S_1^k\, R_0^p\, R_1^q\, X^s$ for some $h, k \in \mZ_n$ and $p, q, s \in \{0,1\}$ we can just apply the relations in (a) to put any word $W$ in our generators into this form.  First move any $X$ to the end of the word, to write $W = W_1 X^s$ with $s \in \{0,1\}$. Then move any $R_0$ or $R_1$ to the end of $W_1$, then rearrange them to write $W_1 = W_2 R_0^p R_1^q$.  Finally $W_2$ can be rearranged to have the form $W_2 = S_0^h S_1^k$.

Now consider the effect of
$T= S_0^h\, S_1^k\, R_0^p\, R_1^q\, X^s$
on the single point $(0.1, 0.2) \in \mR_n \times \mR_n$.
If $s=0$ then
$T(0.1, 0.2) = ((-1)^p\, 0.1 + h, (-1)^q\, 0.2+k)$
and if $s=1$ then
$T(0.1,0.2) = ((-1)^p\, 0.2 + h, (-1)^q\, 0.1 + k)$.
It is easy to determine the values of $s$ then $p, q, h, k$ from $T(0.1, 0.2)$ and therefore $T$ cannot have a different representation.
(This is one place where it is helpful to consider permutations of $\pr_n$, not just $\pz_n$.)

Since there are $n$ choices for each of $h$ and $k$ and $2$ choices for each of $p$, $q$, and $s$, we get $|\Wgd| = 8n^2$.
\end{proof}

\subsection{Dipole coset averages}\label{ca-pf-dipole}

We are going to consider the action of $\Wgd = \la S_0, S_1, R_0, R_1, X \ra \le \sym(\pr_n)$ and some of its subgroups on $\cP_n$, and the corresponding actions on $\sym(\mZ_n)$ and $\cD_n$.  We know that $|\Wgd| = 8n^2$ even if $n=1$ or $2$.  When $n=1$ or $2$ the action of $\Wgd$ on $\cD_n$ is not faithful (two different elements of $\Wgd$ may act in the same way), but this does not matter for counting arguments using Burnside's Lemma.

For determining fixed points we use the action of $\Wgd$ or its subgroups on $\symz_n$.
From Theorem \ref{dsymmetries} and results in the last subsection, we know that the action of an element of $\Wgd$ on $\pi \in \symz_n$ can be described as
\begin{tightequation}
  S_0^h\, S_1^k\, R_0^p\, R_1^q\, X^s (\pi)
	= \si^k \rho^q\, \pi^{(-1)^s} \rho^p \si^{-h}.
\end{tightequation}

\noindent
The action of $\Wgd$ itself will allow us to count dipolar cogs, problem (D6) in our list above.

The subgroups of $\Wgd$ that we consider will always contain $S_0$ and $S_1$, and hence will always have $\Csd = \la S_0, S_1 \ra$ as a subgroup.  From Theorem \ref{dsymmetries} we know that $\Csd = \{S_0^h S_1^k \;|\; h, k \in \mZ_n\}$, with $n^2$ elements.  The action of $\Csd$ on $\cD_n$ or $\Sz_n$ will allow us to count vertex-labeled oriented dipole embeddings, problem (D1) in our list above.

Any group of symmetries containing $\Csd$ can be considered as a union of right cosets $\Csd T$ of $\Csd$.
In $\Wgd$ there are eight right cosets $\Csd T$, for $T \in  \{ R_0^p R_1^q X^s \;|\; p, q, s \in \{0,1\}\} = \{I, X, R, RX, R_0, R_1, R_0 X, R_1 X\}$, where $R = R_0 R_1$.
Therefore, in applying Burnside's Lemma we can use expressions giving the average number of fixed points for $\Csd T$,
$$\Adx(T) = \frac{1}{n^2} \sum_{\ga \in \Csd T} |\fix(\ga)|
	= \frac{1}{n^2} \sum_{h \in \mZ_n} \sum_{k \in \mZ_n}
		|\fix(S_0^h S_1^k T)|
$$
where we are considering the action on $\Sz_n$.
In this subsection we compute all eight corresponding values $\Adx(T)$. We will see that there are five distinct values.

For all of our other dipole counting problems, the relevant group $\Ga$ satisfies $\Csd \le \Ga \le \Wgd$, and the cosets of $\Csd$ in $\Ga$ will be a subset of the cosets of $\Csd$ in $\Wgd$, so we will be able to immediately solve those problems as well.

We can reduce one counting problem to another if two cosets are related by conjugacy, as follows.

\begin{lemma}\label{conjugatecosets}
Suppose $T_1, T_2, U \in \Wgd$ and $U\iv \Csd T_1 U = \Csd T_2$.  Then $\Adx(T_1)=\Adx(T_2)$.
\end{lemma}

\begin{proof}
The map $\ga \mapsto \ga' = U\iv \ga U$ is a bijection from $\Csd T_1$ to $\Csd T_2$.  Moreover, $\pi \in \fix(\ga)$ if and only if $U\iv \pi \in \fix(\ga')$.  So there is a bijection between $\Csd T_1$ and $\Csd T_2$ that preserves the number of fixed points of each element, and hence $\Adx(T_1) = \Adx(T_2)$.
\end{proof}

\begin{computing}{$\Adx(I)$}
To compute $\Adx(I)$ we consider fixed points of elements of $\Csd$. Suppose there is $\pi \in \fix(S_0^h S_1^k)$. Then
$\pi = \si^k \pi \si^{-h}$, so
$\pi \si^h = \si^k \pi$, and thus
$\pi(i+h) = \pi(i) + k$ for all $i \in \mZ_n$.
Hence, by induction,
\begin{equation}\label{eq:dI}
  \pi(i+th) = \pi(i)+tk \quad\text{for all $t \ge 1$ and $i \in \mZ_n$.}
\end{equation}

Let $g = (h,n)$ and $d = n/g$, and $j = (k,n)$ and $e = n/j$.  Then $d$ and $e$ are the smallest positive integers such that $n \dv dh$ and $n \dv ek$, respectively.  If $d < e$ we have $\pi(i) = \pi(i+dh) = \pi(i)+dk \ne \pi(i)$, and if $e < d$ we have $\pi(i) \ne \pi(i+eh) = \pi(i)+ek = \pi(i)$, both of which are contradictions.  Hence $d=e$ and so $g=j$, i.e., $(h,n)=(k,n)$; otherwise no such $\pi$ exists.

Now if $(h,n)=(k,n)=g$ then we can write $h=ag$, $k=bg$, and $n=dg$, where $(a,d)=(b,d)=1$.  For a given $g$, there are $\eul(d)$ possible values of $a$ (and hence of $h$), and the same number of values of $b$ (and hence of $k$).

Consider one of the $\eul(d)^2$ pairs $(h, k)$ corresponding to a given $g$ and $d$.  Now $h\mZ_n = k\mZ_n = g\mZ_n$, which has $d = n/g$ elements.  
By equation \eqref{eq:dI}, once we determine $\pi(i)$ we determine $\pi$ for every element of the coset $i + h\mZ_n = i + g\mZ_n$, and those values exhaust all the elements of the coset $\pi(i) + k\mZ_n = \pi(i) + g\mZ_n$.
Thus, $\pi \in \fix(S_0^h S_1^k)$ is determined by its values $\pi(0), \pi(1), \dots, \pi(g-1)$, and each of these values lies in a distinct coset $j+g\mZ_n$ for $0 \le j \le g-1$, each of which has $d$ elements.  We may therefore assign the cosets to $\pi(0), \pi(1), \dots, \pi(g-1)$ in $g!$ ways, and then pick one of $d$ values in each coset, giving $g!\, d^g$ choices of $\pi$.

Therefore,
$$\Adx(I)
	= \frac{1}{n^2} \sum_{h \in \mZ_n} \sum_{k \in \mZ_n}
		|\fix(S_0^h S_1^k)|
   = \frac{1}{n^2} \sum_{(d,g)\,:\, dg = n} \mskip-10mu \eul(d)^2\, g!\, d^g
   = \frac{1}{n} \sum_{(d,g)\,:\, dg = n} \mskip-10mu \eul(d)^2\, (g-1)!\, d^{g-1}.
$$
\end{computing}

\begin{computing}{$\Adx(X)$}
Suppose $\ga = S_0^h S_1^k X \in \Csd X$ and $\pi \in \fix(\ga)$.  Then $\pi = \si^k \pi\iv \si^{-h}$, from which $(\pi \si^h)^2 = \si^{h+k}$. Thus, the number of fixed points $\pi$ is the same as the number of $\tau = \pi \si^h$ that satisfy $\tau^2 = \si^{h+k} = \si^\ell$, where we let $\ell = h+k$.  Each value of $\ell$ occurs for $n$ pairs $(h,k)$.  For a given $\ell$, there are $g=(\ell,n)$ cycles of $\si^\ell$, each with length $d = n/g$.
Therefore, by Lemma \ref{squareroots}, there are $\sr(d,g)$ possible $\tau$ and hence $\sr(d,g)$ possible $\pi$.  Writing $\ell = cg$ and $n = dg$, we see that $(c,d) = 1$, so there are $\eul(d)$ possible values of $c$, and hence of $\ell$, for a given $g$ and $d$.   Thus, applying Lemma \ref{squareroots},
\begin{align*}
\Adx(X) &= \frac{1}{n^2} \sum_{h \in \mZ_n} \sum_{k \in \mZ_n}
		|\fix(S_0^h S_1^k X)|
	= \frac{1}{n^2} \sum_{(d,g)\,:\, dg =n}
		n\, \eul(d)\, \sr(d,g) \\
   &= \frac{1}{n} \sum_{\substack{(d,g)\,:\, dg = n \\ \text{$d, g$
even}}}
		\eul(d)\, \ma(g, g/2)\, d^{g/2}
	+ \frac{1}{n} \sum_{\substack{(d,g)\,:\, dg = n \\ \text{$d$ odd}}}
		\eul(d) \sum_{j=0}^{\lfloor g/2 \rfloor} \ma(g, j)\, d^j.
\end{align*}
Note that the first term here is nonzero only if $n \equiv 0 \pmod4$.

\end{computing}

\begin{computing}{$\Adx(R)$}
Suppose $\ga = S_0^h S_1^k R = S_0^h S_1^h R_0 R_1$ and $\pi \in \fix(\ga)$.  Then 
$\pi = \si^k \rho \pi \rho \si^{-h}$, from which
$\pi(i) = k-\pi(-(i-h)) = k - \pi(h-i)$ and hence
$\pi(i) + \pi(h-i) = k$ for all $i \in \mZ_n$.
Our analysis of this equation will depend on whether $n$ is odd or even.

Suppose first that $n$ is odd.  Then the sets $\{i, h-i\}$ partition $\mZ_n$ into $(n-1)/2$ pairs and one singleton $\{h/2\}$ (since $n$ is odd, $h/2$ is well-defined for all $h \in \mZ_n$).  We must have $\pi(h/2) + \pi(h-h/2) = k$, i.e., $2\pi(h/2) = k$, so $\pi(h/2) = k/2$ is determined.  For $i \ne h/2$, the pair $\pi(\{i, h-i\}) = \{\pi(i), \pi(h-i)\}$ is a pair $\{j, k-j\}$ with $j \ne k/2$.  Since $\pi$ maps $(n-1)/2$ pairs to another $(n-1)/2$ pairs, by Observation \ref{matchingpairs}, the number of possible $\pi$ is $((n-1)/2)!\, 2^{(n-1)/2}$.
This is constant for all $h$ and $k$.  Therefore, we have
$$\Adx(R) = \left(\frac{n-1}{2}\right)!\; 2^{(n-1)/2}
	\quad\text{for odd $n$}.$$

Suppose now that $n$ is even.  Then the way in which the sets $\{i,h-i\}$ partition $\mZ_n$ depends on whether $h$ is odd or even.  If $h$ is odd, the partition has $n/2$ pairs.  If $h$ is even, define $h/2$ by treating $h \in [0,n-1]$ as a real number rather than an element of $\mZ_n$; then there are two singletons $\{h/2\}$ and $\{h/2+n/2\}$, and $n/2-1$ pairs.  Similarly, the sets $\{j, k-j\}$ partition $\mZ_n$ in ways that depend on whether $k$ is odd or even.  Since each set $\{i, h-i\}$ must map to a set $\{j, k-j\}$ of the same size, either both $h$ and $k$ are odd, or both $h$ and $k$ are even. If both are odd then $\pi$ maps $n/2$ pairs to another $n/2$ pairs, so there are $(n/2)!\, 2^{n/2}$ possible $\pi$ for these $n^2/4$ choices of $(h,k)$.
If both $h$ and $k$ are even there are $2$ ways to match up the singletons $\{h/2\}$ and $\{h/2+n/2\}$ with the singletons $\{k/2\} and \{k/2+n/2\}$, and $(n/2-1)!\, 2^{n/2-1}$ ways for $\pi$ to preserve the pairings of the remaining elements.
So there are $(n/2-1)!\, 2^{n/2}$ possible $\pi$ for these $n^2/4$ choices of $(h,k)$.  Thus,
\begin{align*}
\Adx(R) &= \frac{1}{n^2} \left(
	\frac{n^2}{4} \left(\frac{n}{2}\right)!\, 2^{n/2} +
	\frac{n^2}{4} \left(\frac{n}{2}-1\right)!\, 2^{n/2}
	\right) \\
 &= (n+2) \left(\frac{n}{2}-1\right)!\; 2^{n/2-3}
	\quad\text{for even $n$}.
\end{align*}
\end{computing}

\begin{computing}{$\Adx(RX)$}
By Theorem \ref{dsymmetries} we can see that $RX = R_0 R_1 X = R_0 X R_0$, and also that $\Csd R_0 = R_0 \Csd$.  Therefore, $\Csd R X = \Csd R_0 X R_0 = R_0 \Csd X R_0 = R_0\iv \Csd X R_0$.  Hence, by Lemma \ref{conjugatecosets} we have $$\Adx(RX) = \Adx(X).$$ \end{computing}

Since $R_1$ affects the output of a permutation, and $R_0$ affects the input, the effect of $R_1$ is easier to analyze than the effect of $R_0$.  Therefore, we will consider $R_1$ before $R_0$, and $R_1 X$ before $R_0 X$.

\begin{computing}{$\Adx(R_1)$}
Suppose $\ga = S_0^h S_1^k R_1$ and $\pi \in \fix(\ga)$.  Then
$\pi = \si^k \rho \pi \si^{-h}$, so
$\pi(i) =  k-\pi(i-h)$ and hence
$\pi(i) + \pi(i-h) = k$ for all $i \in \mZ_n$.
Then $\pi(i-h) + \pi(i-2h) = k$ for all $i$ so $\pi(i-2h) = \pi(i)$ for all $i$.  This can only happen if $2h=0$ which means $h=0$, or $n$ is even and $h = n/2$.

If $h=0$ then the only way we can have $\pi(i) + \pi(i-h) = 2\pi(i) = k$ for all $i$ is if $n=1$ or $2$ and $k=0$.  There are one such $\pi$ for $n=1$ and two such $\pi$ for $n=2$.

If $n$ is even and $h = n/2$ then $\pi$ must map each of the $n/2$ pairs $\{i, i+n/2\}$ to a pair $\{j, k-j\}$. For all sets $\{j, k-j\}$ to be pairs, $k$ must be odd.  For each of the $n/2$ odd values of $k$, Observation \ref{matchingpairs} tells us there are $(n/2)!\, 2^{n/2}$ possible $\pi$.
Thus, taking into account the special cases when $n = 1$ or $2$,
\begin{align*}
 \Adx(R_1) 
   &= \begin{cases}
     1 & \text{if $n=1$ or $2$,}\\
     0 & \text{if $n \ge 3$ is odd,}\\
 \displaystyle \frac{1}{n} \left( \frac{n}{2} \right)! \; 2^{n/2-1}
	& \text{if $n\ge 4$ is even.}
   \end{cases}
\end{align*}
\end{computing}

\begin{computing}{$\Adx(R_0)$}
By Theorem \ref{dsymmetries} we can see that $R_0 = X R_1 X$, and that $\Csd X = X \Csd$.  Therefore, $\Csd R_0 = \Csd X R_1 X = X \Csd R_1 X = X\iv \Csd R_1 X$.  Hence, by Lemma \ref{conjugatecosets} we have $$\Adx(R_0) = \Adx(R_1) .$$
\end{computing}

\begin{computing}{$\Adx(R_1 X)$}
Suppose $\ga = S_0^h S_1^k R_1 X$ and $\pi \in \fix(\ga)$.  Then $\pi = \si^k \rho \pi\iv \si^{-h}$, from which $(\si^h \pi)^2 = \si^{h+k} \rho$. Thus, the number of fixed points $\pi$ is the same as the number of $\tau = \si^h \pi$ that satisfy $\tau^2 = \si^{h+k} \rho = \si^\ell \rho = \al$, where we let $\ell = h+k$.  Each value of $\ell$ occurs for $n$ pairs $(h,k)$.  Now $\al(i) = \si^\ell \rho(i) = \ell-i$, so $\al$ is always an involution.

If $n$ is odd then for all $n$ values of $\ell$, then the number of $2$-cycles of $\al$ is $(n-1)/2$ and $\al$ has one $1$-cycle, namely $(\ell/2)$.  Applying Lemma \ref{squareroots} we get
\begin{align*}
 \Adx(R_1 X) &=  \frac{1}{n^2} n^2\, \sr(1, 1) \sr(2, (n-1)/2)
		= \sr(2, (n-1)/2) \\[4pt]
   &= \begin{cases}
	0 & \text{if $(n-1)/2$ is odd, i.e., $n \equiv 3 \pmod4$,}\\
	\ma((n-1)/2, (n-1)/4)\, 2^{(n-1)/4}
		& \text{if $(n-1)/2$ is even, i.e.,
			$n \equiv 1 \pmod4$.}
    \end{cases}
\end{align*}

If $n$ is even then the number of $2$-cycles of $\al$ is $n/2$  for the $n/2$ odd values of $\ell$.  For the $n/2$ even values of $\ell$, the number of $2$-cycles of $\al$ is $n/2-1$, and there are also two $1$-cycles, namely $(\ell/2)$ and $(\ell/2+n/2)$.  Applying Lemma \ref{squareroots} we get
\begin{align*}
 \Adx(R_1 X) &= \frac{1}{n^2} \left(
	  \frac{n^2}{2} \sr(2, n/2)
	  + \frac{n^2}{2} \sr(1, 2) \sr(2, n/2-1)
	\right) \\[4pt]
   &= \begin{cases}
	\frac{1}{2} \ma(n/2, n/4) 2^{n/4}
		& \text{if $n/2$ is even, i.e., $n \equiv 0 \pmod4$,}\\
	\frac{1}{2} 2 \ma(n/2-1, \frac{1}{2}(n/2-1)) 2^{\frac{1}{2}(n/2-1)}
		& \text{if $n/2-1$ is even, i.e., $n \equiv 2 \pmod4$,}\\
     \end{cases} \\[4pt]
   &= \begin{cases}
	\ma(n/2, n/4)\, 2^{n/4-1}
		& \text{if $n \equiv 0 \pmod4$,}\\
	\ma((n-2)/2, (n-2)/4)\, 2^{(n-2)/4}
		& \text{if $n \equiv 2 \pmod4$,}\\
     \end{cases}
\end{align*}

\end{computing}

\begin{computing}{$\Adx(R_0 X)$}
By Theorem \ref{dsymmetries} we can see that $R_0 X = X R_1$, and that $\Csd X = X \Csd$.  Therefore, $\Csd R_0 X = \Csd X R_1 = X \Csd R_1 = X (\Csd R_1 X) X\iv$.
Hence, by Lemma \ref{conjugatecosets} we have
$$\Adx(R_0 X) = \Adx(R_1 X) .$$
\end{computing}

\subsection{Counting dipole embeddings and related objects}\label{bo-pf-dipole}

In this subsection we prove counting formulas (D1)--(D8) from Subsection \ref{res-dipole}.
We consider $\cD_n$, or equivalently $\Sz_n$, under the action of various groups $\Ga$ with $\Csd \le \Ga \le \Wgd$.  Each such group can be written as a union of cosets of $\Csd$, i.e., as $\Ga = \Csd T_1 \wcup \Csd T_2 \wcup \dots \wcup \Csd T_k$ for some $T_i \in \la R_0, R_1, X\ra = \{R_0^p R_1^q X^s \;|\; p, q, s \in \{0,1\}\}$.  Therefore, by Burnside's Lemma the number of equivalence classes is
\begin{equation}\label{usecosetav}
\frac{1}{|\Ga|} \sum_{\ga \in \Ga} |\fix(\ga)|
   = \frac{1}{kn^2} \sum_{i = 1}^k \sum_{\ga \in \Csd T_i} |\fix(\ga)|
   = \frac{1}{k} \sum_{i = 1}^k \Adx(T_i) .
\end{equation}

\begin{counting}{(D1) vertex-labeled oriented dipole embeddings}
Given a vertex-labeled oriented dipole embedding $\Phi$, we can transform it into a labeled dipole by choosing a half-edge incident with vertex $0$ to label $0$, and then labeling the other half-edges incident with vertex $0$ in ascending clockwise order, and similarly for vertex $1$.  However, our choices of which half-edge incident with vertex $0$ to label $0$, and which half-edge incident with vertex $1$ to label $0$, are arbitrary.  So other labeled dipoles for $\Phi$ can be obtained by applying arbitrary cyclic shifts $S_0^h$ and $S_1^k$.

Thus, a vertex-labeled oriented dipole embedding may be regarded as an equivalence class of labeled dipoles under the action of $\Csd = \la S_0, S_1 \ra = \Csd I$, and so the number of equivalence classes is just $\Adx(I)$.
\end{counting}

\begin{counting}{(D2) oriented dipole embeddings}
To count oriented dipole embeddings where the vertices are unlabeled, we include vertex swaps in our allowed symmetries, so the group of symmetries is $\la S_0, S_1, X \ra$, which from Theorem \ref{dsymmetries} we know is $\{S_0^h S_1^k X^s \;|\; h, k \in \mZ_n, s \in \{0,1\}\} = \Csd \wcup \Csd X$.
The number of equivalence classes is therefore
$\rc2 (\Adx(I)+\Adx(X))$.
\end{counting}

\begin{counting}{(D3) vertex-labeled orientable dipole embeddings}
To count orientable dipole embeddings we must allow for the orientation of the surface being reversed.  This means that when we choose the two half-edges incident with vertices $0$ and $1$ to give label $0$, we then label the other half-edges in the reverse order, at both vertices.  This means we are applying the action of $R = R_0R_1$.
The group of symmetries in the vertex-labeled case is therefore $\la S_0, S_1, R\ra$, which by applying Theorem \ref{dsymmetries} can be written as
$\{S_0^h S_1^k R^r \;|\; h, k \in \mZ_n, r \in \{0,1\}\}
	= \Csd \wcup \Csd R$.
The number of equivalence classes is therefore
$\rc2 (\Adx(I)+\Adx(R))$.
\end{counting}

\begin{counting}{(D4) orientable dipole embeddings}
To count orientable dipole embeddings where the vertices are unlabeled, we add $X$ to the group of symmetries, giving 
$\la S_0, S_1, R, X\ra$. By Theorem \ref{dsymmetries} this can be
written as
$\{S_0^h S_1^k R^r X^s \;|\; h, k \in \mZ_n, r, s \in \{0,1\}\}
	= \Csd \wcup \Csd R \wcup \Csd X \wcup \Csd RX$.  
The number of equivalence classes is therefore
$\rc4 (\Adx(I) + \Adx(R) + \Adx(X) + \Adx(RX))$.
Since $\Adx(RX) = \Adx(X)$, this simplifies to 
$\rc4 (\Adx(I) + \Adx(R) + 2\Adx(X))$.
\end{counting}

\begin{counting}{(D5) vertex-labeled dipolar cogs}
To count cogs we must allow the cyclic ordering at each vertex to be reversed independently.  This means that we can apply the actions of both  $R_0$ and $R_1$.  The group of symmetries in the vertex-labeled case is therefore
$\la S_0, S_1, R_0, R_1 \ra$,
which by Theorem \ref{dsymmetries} can be written as
$\{S_0^h S_1^k R_0^p R_1^q \;|\; h, k \in \mZ_n, p, q \in \{0,1\}\}
	= \Csd \wcup \Csd R_0 \wcup \Csd R_1 \wcup \Csd R_0R_1$.  
The number of equivalence classes is therefore
$\rc4 (\Adx(I)+\Adx(R_0)+\Adx(R_1)+\Adx(R_0R_1))$.
Since $\Adx(R_0)=\Adx(R_1)$ and $R_0R_1 = R$, this simplifies to
$\rc4 (\Adx(I) + \Adx(R) + 2\Adx(R_1))$.
\end{counting}

\begin{counting}{(D6) dipolar cogs}
Again, when the vertices are unlabeled we add $X$ to the group of symmetries, which is therefore
$\Wgd = \la S_0, S_1, R_0, R_1, X\ra$.
By Theorem \ref{dsymmetries} this can be written as
$\{S_0^h S_1^k R_0^p R_1^q X^s \;|\; h, k \in \mZ_n, p, q, s \in
		\{0,1\}\}
	= \Csd \wcup \Csd R_0 \wcup \Csd R_1 \wcup \Csd R_0R_1 \wcup
	 \Csd X \wcup \Csd R_0X \wcup \Csd R_1X \wcup \Csd R_0R_1X $.
The number of equivalence classes is therefore
\begin{center}
$\rc8 (\Adx(I)+\Adx(R_0)+\Adx(R_1)+\Adx(R_0R_1)+
	\Adx(X)+\Adx(R_0X)+\Adx(R_1X)+\Adx(R_0R_1X))$.
\end{center}
Since $\Adx(R_0)=\Adx(R_1)$, $\Adx(R_0X) = \Adx(R_1X)$, $R_0R_1 = R$ and $\Adx(R_0R_1 X) = \Adx(RX) = \Adx(X)$, this simplifies to
$\rc8 (\Adx(I) + \Adx(R) + 2\Adx(R_1) + 2\Adx(X) + 2\Adx(R_1X))$.
\end{counting}

\begin{counting}{(D7) equivalence classes of permutations under cyclic shifts and reversal of input variables only (or output variables only)}
There are two subgroups of $\Wgd$ that do not have very natural interpretations in terms of dipole embeddings, but can be considered as groups of symmetries of $\Sz_n$.  These are where we allow cyclic shifts of both input and output variables, and reversal of just one set of variables (input or output, but not both).
These give the groups $\Ga_7 = \la S_0, S_1, R_1\ra$ and $\Ga_7' = \la S_0, S_1, R_0\ra$.  Since $\Ga_7' = X\iv \Ga_7 X$, the number of equivalence classes for both groups will be the same, by applying Lemma \ref{conjugatecosets}.  So we just consider $\Ga_7 = \la S_0, S_1, R_1 \ra$, which by Theorem \ref{dsymmetries} can be written as
$\{S_0^h S_1^k R_1^q \;|\; h, k \in \mZ_n, q \in \{0,1\}\}
	= \Csd \wcup \Csd R_1$.
The number of equivalence classes is therefore
$\rc2 (\Adx(I)+\Adx(R_1))$.
\end{counting}

\begin{counting}{(D8) equivalence classes of permutation matrices under cyclic shifts of the row set, cyclic shifts of the column set, and rotations by multiples of $90^\circ$}
For permutation matrices, a cyclic shift of the row set corresponds to $S_0$, a cyclic shift of the column set corresponds to $S_1$, and a rotation by $90^\circ$ clockwise corresponds to $V_1 X = S_1^{n-1} R_1 X$ (a transposition, $X$, followed by reversal of the column set, $V_1 = S_1^{n-1} R_1$).
Thus, the group of symmetries is $\la S_0, S_1, R_1 X\ra$.  By Theorem \ref{dsymmetries} we see that $(R_1X)^2 = R$, $(R_1 X)^3 = R_0 X$, and $(R_1 X)^4 = I$.  So the group can be written as $\Csd \wcup \Csd R \wcup \Csd R_0 X \wcup \Csd R_1 X$.  The number of equivalence classes is therefore
$\rc4( \Adx(I) + \Adx(R) + \Adx(R_0 X) + \Adx(R_1 X) )$.
Since $\Adx(R_0 X) = \Adx(R_1 X)$, this simplifies to
$\rc4( \Adx(I) + \Adx(R) + 2 \Adx(R_1 X) )$.
\end{counting}

There are ten groups $\Ga$ that satisfy $\Csd \le \Ga \le \Wgd$.  These correspond to the subgroups of the quotient group $\Wgd/\Csd$, which is an $8$-element dihedral group.  Thinking of $\Wgd/\Csd$ as symmetries of a square gives a natural correspondence with operations on permutation matrices.

We have counted equivalence classes of $\cD_n$ (or $\cP_n$ or $\Sz_n$) under the action of $\Ga$ for eight of these groups in (D1)--(D8).  Actually, as we noted above, (D7) handles two of these groups that are conjugate in $\Wgd$.  Item (D2) also handles two of these groups, because the groups  $\Ga_2 = \la S_0, S_1, X \ra$ and $\Ga_2' = \la S_0, S_1, RX \ra$ are conjugate in $\Wgd$, with $\Ga_2' = R_0\iv \Ga_2 R_0$.
Thus, we have covered all ten groups.
Since $X R_0 X = R_1$, the groups $\la S_0, S_1, R_0, X\ra$ and $\la S_0, S_1, R_1, X\ra$ are just the full group $\Wgd = \la S_0, S_1, R_0, R_1, X\ra$, and therefore we do not need to consider separate counting questions involving these groups.

\begin{observation}
The above results tell us that all of $\Adx(I)$, $\Adx(R)$, $\Adx(R_1)$, $\Adx(X)$, and $\Adx(R_1X)$ are integers, as follows.  By (D1) we know that $\Adx(I)$ is the number of equivalence classes under $\Csd$, so it is an integer.
By (D2), (D3), and (D7) we also know that $\Adx(I)+\Adx(X)$, $\Adx(I)+\Adx(R)$, and $\Adx(I)+\Adx(R_1)$ are even integers.  Therefore, $\Adx(X)$, $\Adx(R)$, and $\Adx(R_1)$ are integers.
Finally, we know by (D6) that $\Adx(I)+\Adx(R)+2\Adx(R_1)+2\Adx(X)+ 2\Adx(R_1) + 2\Adx(R_1X)$ is divisible by $8$, so is even, and all of $\Adx(I)+\Adx(R)$, $2\Adx(R_1)$, and $2\Adx(X)$ are even.  Thus, $2\Adx(R_1X)$ is even and $\Adx(R_1 X)$ is an integer.
\end{observation}

\section{Proofs for bouquet formulas}\label{pf-bouquet}

\subsection{Colored labeled bouquets and symmetry operations}\label{so-pf-bouquet}

In this section we prove counting formulas (B1)--(B4) from Subsection \ref{res-bouquet} regarding embeddings of colored bouquets. Recall that a bouquet $B_n$ has one vertex and $n$ loops.  Again we will think of each edge as consisting of two half-edges.  By `colored' we mean that each edge receives an arbitrary color from a set of $k$ colors.
By using results in the case $k=2$ we are able to count nonorientable embeddings of bouquets, where previous counting results for embeddings of bouquets have only considered orientable embeddings.  Our results will be proved by elementary techniques based on groups acting on a set of objects that we will call colored labeled bouquets.

A \emph{colored labeled bouquet} is a bouquet $B$ where the half-edges received distinct labels from $\Ztn = \{0, 1, 2, \dots, 2n-1\}$, where $n = |E(B)|$, and where each edge receives a color from $\mZ_k = \{0, 1, 2, \dots, k-1\}$.  We let $\cB_{n,k}$ denote the set of $n$-edge $k$-colored bouquets.

A colored labeled bouquet $B$ is completely described by a  perfect matching $M$ in the complete graph $K(\Ztn)$ with vertex set $\Ztn$, plus a coloring function $\col : M \to \mZ_k$. For each $e \in E(B)$, the perfect matching $M$ contains an edge $\set{a_e, b_e}$ (describing an edge as an unordered pair of vertices), where $a_e$ and $b_e$ are the labels of the two half-edges of $e$, and $\col$ assigns the color of $e \in E(B)$ to $\set{a_e, b_e} \in M$.  There is a one-to-one correspondence between $\cB_{n,k}$ and $\cM_{n,k}$, the set of $k$-colored perfect matchings $(M, \col)$ in $K(\Ztn)$. Perfect matchings on a cyclically ordered set such as $\Ztn$ are often known as \emph{chord diagrams}, so our results may also be interpreted as results for colored chord diagrams.  See Figure \ref{fig:bouquetchord}.

Each oriented or orientable embedding of a colored bouquet can be turned into a colored labeled bouquet (or $k$-colored perfect matching) in a natural way.  We refer the reader back to Figure \ref{fig:Interleaving} for an example. But the colored labeled bouquet is not in general unique, and we wish to characterize the possible bouquets by equivalence under certain symmetry operations.  The symmetries always include cyclically shifting the labels of the half-edges, but may also include reversing the labels.

Our symmetries can be defined in terms of their effect on elements of $\Ztn$, but to avoid treating $n=1$ and $2$ as special cases, and to simplify some proofs, we define them more generally as permutations of $\mR_{2n} = \mR/2n\mR$, which contains $\Ztn$ as a subgroup.   We define $S, R \in \symr_{2n}$ as follows:

 \listitem $S(a) = a+1$;\quad
     and\quad $R(a) = -a$.

\noindent
We let $\Wgb = \la S, R \ra \subseteq \symz_{2n}$.

If $T = S$ or $R$ then it is clear that applying $T$, or $T\iv$, to all the labels in a colored labeled bouquet, without changing any colors, produces a new colored labeled bouquet, so that $T$ and $T\iv$ permute $\cB_{n,k}$.
Correspondingly, we may apply $T$, or $T\iv$, to the vertices of $K(\Ztn)$, and define $T(M) = \{ \set{T(a),T(b)} \st \set{a,b} \in M\}$, with the edge $T(\set{a,b}) = \set{T(a), T(b)}$ receiving the color assigned to $\set{a,b}$.  In this way $T$ and $T\iv$ permute $\cM_{n,k}$. We can then extend this to all $T \in \Wgb$, so we have an action of $\Wgb$ on $\cB_{n,k}$ and an equivalent action on $\cM_{n,k}$.

The basic properties of $S$, $R$ and $\Wgb$ are as follows.

\begin{thm}\label{bsymmetries}
Let $n$ be a positive integer, and consider $S, R \in \sym(\mR_{2n})$ as defined above.  Let $I$ be the identity of $\sym(\mR_{2n})$.

\smallskip\noindent
(a) Then (composing functions right to left) $S^{2n} = R^2 = I$ and $RS = S\iv R$.

\smallskip\noindent
(b) Every element of $\Wgb = \la S, R \ra$ can be written uniquely as $S^h R^r$ where $h \in \Ztn$ and $r \in \{0,1\}$.  Hence, $|\Wgb| = 4n$.
\end{thm}

We omit the proof, which is straightforward.

\subsection{Bouquet coset averages}\label{ca-pf-bouquet}

We are going to consider the action of $\Wgb = \la S, R \ra \le \symr_{2n}$ (a dihedral group) and its subgroup $\Csb = \la S \ra$ (a cyclic group) on $\cB_{n,k}$ and $\cM_{n,k}$.  When $n=1$ or $2$ the action is not faithful (for example, if $n=2$ then $S^2$ acts in the same way as $I$) but this does not matter.

For determining fixed points for Burnside's Lemma we use the action of $\Wgb$ or $\Csb$ on $\cM_{n,k}$. As in Section \ref{pf-dipole}, we will consider average numbers of fixed points for the right cosets of $\Csb$, which are just $\Csb$ itself and $\Csb R$.  Each coset has size $|\Csb| = 2n$ so for $T = I$ or $R$ we define
$$\Abx(T) = \frac{1}{2n} \sum_{\ga \in \Csb T} |\fix(\ga)|
	= \frac{1}{2n} \sum_{h \in \Ztn}
		|\fix(S^h T)|
$$
where we are considering the action on $\cM_{n,k}$.

Before computing coset averages we consider what a fixed point $(M,\col) \in \cM_{n,k}$ of $T \in \Wgb$ must look like in general.  Suppose $\set{a,b} \in M$.  Then we also have $\set{T(a), T(b)} \in M$ and hence we have $\set{T^t(a), T^t(b)} \in M$ for all integers $t$.  For $t \ge 0$, this follows by induction.  For $t < 0$ it follows because something is a fixed point of $T$ if and only if it is a fixed point of $T\iv$.  The edges $\set{T^t(a), T^t(b)} \in M$ for all integers $t$ form an edge-orbit of $K(\Ztn)$ under the action of $\la T \ra$.  Since $\col$ is also fixed by $T$, all edges in an edge-orbit must have the same color.

The points $T^t(a)$ form a vertex-orbit of the action of $\la T \ra$ on $K(\Ztn)$, as do the points $T^t(b)$.  If the two vertex-orbits are the same, we must have an orbit of even size, because $M$ matches up the vertices in the vertex-orbit.  If the two vertex-orbits are different, they must have the same size, because $M$ matches the vertices of one vertex-orbit to the vertices of the other.

Thus, if $(M, \col) \in \cM_{n,k}$ is a fixed point of $T$, then the edges of $M$ can be partitioned into edge-orbits under the action of $\la T \ra$.  Each edge-orbit either matches the vertices of a single vertex-orbit $O$ of even size, or matches the vertices of one vertex-orbit $O_1$ to the vertices of a paired vertex-orbit $O_2$ of the same size, say $s$. In the latter case, if we fix $a \in O_1$ then there are $s$ possible vertices $b \in O_2$ for which we could have $\set{a,b} \in M$, and once we choose this edge the rest of the edge-orbit is determined. So (assuming no special restrictions apply) there are $s$ possible edge-orbits matching two paired vertex-orbits. All edges of $M$ in the same edge-orbit receive the same color under $\col$.

\begin{computing}{$\Abx(I)$}
To compute $\Abx(I)$ we consider fixed points of elements of $\Csb$, which have the form $S^h$.  Suppose there is $(M, \col) \in \fix(S^h)$.  If we let $g = (h, 2n)$ then $h\Ztn = g\Ztn$.  The number of elements of $g\Ztn$ is $d = 2n/g$.  For a given $g$ there are $\eul(d)$ possible values of $h$.

The vertex-orbits of $K(\Ztn)$ under the action of $\la S^h \ra$ are just the cosets $a + h\Ztn = a + g\Ztn$.  So each edge-orbit must match up vertices in the same coset, or must match one coset to another.

Suppose an edge-orbit matches vertices in the same coset $a+g\Ztn$. Then we have an edge $\set{a,a+ug}$ for some integer $u$, where $ug \ne 0$.  We also have an edge $\set{a+ug, a+2ug}$, so we must have $a = a+2ug$, so that $2ug = 0$.  Since $ug \ne 0$, we must have $ug = n$.  Since $n \in g\Ztn$, the number of elements in $g\Ztn$, which is $d$, must be even.

So suppose $d$ is even, so that $n \in g\Ztn$.  We can construct fixed points of $S^h$ as follows.  Choose $j$ with $0 \le j \le \lfloor g/2 \rfloor$, and pair up $j$ pairs of the $g$ cosets, which may be done in $\ma(g, j)$ ways. 
If $a + g\Ztn$ is paired with $b + g\Ztn$, which both have size $d$, then there are $d$ possible choices of edge-orbit, so there are $d^j$ possible edge-orbits for all the paired cosets.
If $a+g\Ztn$ is not one of the paired cosets, then from above each $a' \in a+g\Ztn$ must be joined by $M$ to $a'+n$, so there is only one possibility for the edges in the unpaired cosets.

There are $j$ edge-orbits between paired cosets, and $g-2j$ edge-orbits inside unpaired cosets.  So there are $g-j$ edge-orbits and hence $k^{g-j}$ choices of $\col$.

The total number of fixed points for a particular value of $h$ with an even value of $d$ is therefore $$\sum_{j=0}^{\lfloor g/2 \rfloor} \ma(g, j) \, d^j \, k^{g-j} .$$

Now suppose $d$ is odd, which means that $g=2n/d$ is even.  From above, we cannot have edges inside a single coset $a+g\Ztn$, so all $g$ cosets are paired up.  This means we just have the case $j=g/2$ from above, so the number of fixed points for a particular value of $h$ with an odd value of $d$ is $\ma(g, g/2) d^{g/2} k^{g/2}$.

Putting these together,
\begin{align*}
\Abx(I) &= \frac{1}{2n} \sum_{h \in \Ztn}
		|\fix(S^h)| \\
 &= \frac{1}{2n} \sum_{\substack{(d,g)\,:\, dg = 2n \\ \text{$d$ odd}}}
		\eul(d)\, \ma(g, g/2)\, d^{g/2}\, k^{g/2}
    + \frac{1}{2n} \sum_{\substack{(d,g)\,:\, dg = 2n \\ \text{$d$
even}}}
	 \eul(d) \sum_{j=0}^{\lfloor g/2 \rfloor} \ma(g, j)\, d^j\, k^{g-j}.
\end{align*}
\end{computing}

\begin{computing}{$\Abx(R)$}
Suppose $\ga = S^h R \in \Csb R$, and $(M, \col) \in \fix(\ga)$. For $i \in \Ztn$, $S^h R (i) = h-i$, so the vertex-orbits of $\la S^h R\ra$ are sets $\set{i, h-i}$.  Such a vertex-orbit may have size $1$ if $i = h-i$, which happens when $h$ is even and $i = h/2$ or $h/2+n$.  Otherwise the vertex-orbits have size $2$.

Suppose $h$ is odd.  Then there are $n$ vertex-orbits, all of size $2$. Thus, all edge-orbits in $M$ either match the two vertices in a single vertex-orbit, or match one vertex-orbit to a paired vertex-orbit.  To construct a fixed point of $S^h R$ we choose $j$ with $0 \le j \le \lfloor n/2 \rfloor$, and pair up $j$ pairs of vertex-orbits of size $2$, which may be done in $\ma(n, j)$ ways.  There are $2$ choices of edge-orbit for each pair of vertex-orbits, for $2^j$ total choices. The unpaired vertex-orbits have one edge of $M$ joining their two vertices, which can only be done in one way.

There are $j$ edge-orbits between paired cosets, and $n-2j$ edge-orbits inside unpaired cosets, so there are $n-j$ edge-orbits and hence $k^{n-j}$ choices of $\col$.

The total number of fixed points for one of the $n$ odd values of $h$ is therefore
$$\sum_{j=0}^{\lfloor n/2 \rfloor} \ma(n,j)\, 2^j\, k^{n-j}.$$

Suppose now that $h$ is even.  Then there are two vertex-orbits $\set{h/2}$ and $\set{h/2+n}$ of size $1$, and $(n-1)$ vertex-orbits of size $2$.  There must be an edge of $M$ joining the two vertex-orbits of size $1$, and there are $k$ choices for the color of this edge. Then we can apply the same analysis as above, replacing $n$ by $n-1$, to the vertex-orbits of size $2$.  So the total number of fixed points for one of the $n$ even values of $h$ is
$$k \sum_{j=0}^{\lfloor (n-1)/2 \rfloor} \ma(n-1,j)\, 2^j\, k^{n-1-j}
   = \sum_{j=0}^{\lfloor (n-1)/2 \rfloor} \ma(n-1,j)\, 2^j\, k^{n-j}.
$$

Putting everything together, we have
\begin{align*}
 \Abx(R) &= \frac{1}{2n} \left(
   n \sum_{j=0}^{\lfloor n/2 \rfloor} \ma(n,j)\, 2^j\, k^{n-j}
   + n \sum_{j=0}^{\lfloor (n-1)/2 \rfloor} \ma(n-1,j)\, 2^j\, k^{n-j}
 \right) \\
    &= \frac{1}{2} \left(
   \sum_{j=0}^{\lfloor n/2 \rfloor} \ma(n,j)\, 2^j\, k^{n-j}
   + \sum_{j=0}^{\lfloor (n-1)/2 \rfloor} \ma(n-1,j)\, 2^j\, k^{n-j}
 \right).
\end{align*}
\end{computing}

\subsection{Counting colored bouquet embeddings and related objects}\label{bo-pf-bouquet}

In this subsection we consider $\cB_{n,k}$, or equivalently $\cM_{n,k}$, under the action of either $\Ga = \Csb$ or $\Ga= \Wgb = \Csb \cup \Csb R$. In a similar way to equation \eqref{usecosetav} we can just take the average over all cosets of $\Csb$ in $\Ga$ of the average number of fixed points in a coset.

\begin{counting}{(B1) oriented embeddings of colored bouquets}
Given an oriented embedding $\Phi$ of a colored bouquet, we can transform it into a colored labeled bouquet by choosing a half-edge to label $0$, and then labeling the other half-edge in ascending clockwise order.  This is illustrated by Figure \ref{fig:Interleaving}.  However, our choice of which half-edge to give label $0$ was arbitrary.  So other colored labeled bouquets for $\Phi$ can be obtained by applying an arbitrary cyclic shift $S^h$.

Thus, an oriented embedding of a colored bouquet may be regarded as an equivalence class of colored labeled bouquets under the action of $\Csb = \la S \ra$, and so the number of equivalence classes is just $\Abx(I)$.

This can also be regarded as the number of colored chord diagrams equivalent under cyclic shifts.
\end{counting}

\begin{counting}{(B2) orientable embeddings of colored bouquets}
To count orientable embeddings of colored bouquets, we must allow for the orientation of the surface being reversed.  This means that when we choose the half-edge to label $0$, we can then label the other half-edges in one of two cyclic orders.  These two labelings are related by the transformation $R$.  So the group of symmetries we must consider is 
$\Wgb = \la S, R \ra = \Csb \cup \Csb R$.
The number of equivalence classes is therefore
$\rc2 (\Abx(I)+\Abx(R))$.
\end{counting}

\begin{counting}{(B3) generic (orientable or nonorientable) embeddings of colored bouquets}
Generic embeddings are described by a rotation scheme together with edge signatures, which describe whether an edge should be considered twisted or not.  For graphs in general this representation is not unique.  However, when a graph has only one vertex the representation is unique, and so a generic embedding can be regarded as a rotation scheme (which may be reversed without changing the embedding) together with edge signatures, which are just a $2$-coloring of the edges.  So generic embeddings of bouquets are in one-to-one correspondence with orientable embeddings of $2$-colored bouquets.  More generally, generic embeddings of $k$-colored bouquets are in one-to-one correspondence with orientable embeddings of $2k$-colored bouquets.  Therefore, the number of generic embeddings of $k$-colored bouquets is
$\rc2 (\Ab(I)(n, 2k)+\Ab(R)(n,2k))$.
\end{counting}

\begin{counting}{(B4) nonorientable embeddings of colored bouquets}
The number of nonorientable embeddings is just the number of generic embeddings minus the number of orientable embeddings.  So the number of nonorientable embeddings of $k$-colored bouquets is
$\text{(B3)}-\text{(B2)} =
	\rc2 (\Ab(I)(n,2k)+\Ab(R)(n,2k)-\Abx(I)-\Abx(R))$.
\end{counting}

\begin{observation}
From (B1) and (B2) we see that $\Abx(I)$ and $\Abx(R)$ are both integers.
\end{observation}

\section{Proofs for directed bouquets}\label{pf-dirbouquet}

\subsection{Colored signed labeled bouquets and symmetry operations}\label{so-pf-dirbouquet}

In this section we prove counting formulas (A1)--(A9) from Subsection \ref{res-dirbouquet} for directed embeddings of colored directed bouquets and related objects. Recall that a directed bouquet $\dg{B}_n$ is a digraph with one vertex and $n$ directed loops.
We will think of each arc (directed edge) as consisting of an outward half-arc and an inward half-arc.  
A directed embedding of a directed bouquet requires that the directions on the half-arcs alternate when going around the cyclic order at the vertex.  By `colored' we mean that each arc receives an arbitrary color from a set of $k$ colors. 

We also consider equivalence classes of digraphs under the operation of reversing the direction of all of the arcs.  We call such an equivalence class an \emph{arc-reversal class}.  If we have a directed embedding of a digraph, reversing all the arcs preserves the fact that we have a directed embedding, so we can also consider arc-reversal classes of directed embeddings.

As with embeddings of bouquets, by using results in the case $k=2$ we are able to count nonorientable directed embeddings of directed bouquets.  Our results will be proved by elementary techiques based on groups acting on a set of objects that we will call colored signed labeled bouquets.

A \emph{colored signed labeled bouquet} is a bouquet $B$ where each of the half-edges receives the following:  a distinct label from $\Ztn = \set{0, 1, 2, \dots,2n-1}$, where $n = |E(B)|$;   either a $+$ or a $-$ sign so that all even half-edges have one sign and all odd half-edges have the opposite sign and so that every edge has both an even half-edge and an odd half-edge; and  a color from $\mZ_k = \set{0, 1, 2, \dots, k-1}$.  The signs indicate how to convert $B$ into a directed bouquet, namely by directing each edge from its positive half-edge to its negative half-edge. Requiring the even half-edges to have one sign and the odd half-edges to have the other assures the alternation of signs required for a directed embedding of the bouquet.  We let $\cA_{n,k}$ denote the set of $n$-edge $k$-colored signed labeled bouquets.

A colored signed labeled bouquet is completely described by a triple $(N, \col, \ep)$.  Here $N$ is a perfect matching in the complete bipartite graph $K(2\Ztn, 1+2\Ztn)$, whose vertices are partitioned into the set $2\Ztn$ of even numbers and the set $1+2\Ztn$ of odd numbers.  The function $\col : N \to \mZ_k$ assigns one of $k$ colors to each edge in $N$.  And finally $\ep \in \{+,-\}$ is the sign assigned to all the even vertices (elements of $2\Ztn$), while $-\ep$ is the sign
assigned to all the odd vertices (elements of $1+2\Ztn$). The vertices of the graph correspond to the half-edges in the bouquet as described in detail in Subsection \ref{so-pf-bouquet}, and the edges of $N$ correspond to the edges of the bouquet.  We let $\cN_{n,k}$ denote the set of triples $(N, \col, \ep)$.

We will use the same symmetries as in Subsection \ref{so-pf-bouquet}, namely $S, R \in \sym(\mR_{2n})$ where $\mR_{2n} = \mR/2n\mR$.
Thus, if $T \in \la S, R \ra$ then we have an action of $T$ on $V(K(2\Ztn, 1+2\Ztn)) = \Ztn \subseteq \mR_{2n}$, which we use to define
$T(N, \col, \ep) = (T(N), \col', \ep')$ where
$T(N) = \{ \set{T(a), T(b)} \;|\ \set{a,b} \in N\}$,
with $\col'(\set{T(a),T(b)}) = \col(\set{a,b})$, and
$\ep'$ is $\ep$ if $T$ preserves the parity of elements of $\Ztn$ and $-\ep$ otherwise.
The reader may think of $T$ as relabeling the vertices without changing edge colors or vertex signs.

However, we have a third basic symmetry $F$, which flips the sign of all vertices.  Clearly $F^2 = I$ and $F$ commutes with elements of $\la S, R \ra$, so our overall group of symmetries is $\Wga = \la S, R \ra \times \la F \ra$. We regard $F$ as acting as the identity on both vertices and edges of $K(2\Ztn, 1+2\Ztn)$.

\subsection{Directed bouquet coset averages}\label{ca-pf-dirbouquet}

We are going to consider the action of  $\Wga = \la S, R \ra \times \la F \ra$ and some of its subgroups on $\cA_{n,k}$ and $\cN_{n,k}$.  For small $n$ the action is not faithful but this does not matter.  All the groups $\Ga$ that we will consider contain $\Csa = \la S \ra$ as a subgroup, so that $\Csa \le \Ga \le \Wga$, and we can write $\Ga$ as a union of cosets of $\Csa$ in $\Wga$.  There are four such cosets $\Csa, \Csa R, \Csa F, \Csa R F$, and for each $T \in \set{I, R, F, RF}$ we will compute the coset average
$$\Aax(T) = \frac{1}{2n} \sum_{\ga \in \Csa T} |\fix(\ga)|
	= \frac{1}{2n} \sum_{h \in \Ztn}
		|\fix(S^h T)|
$$
for the action on $\cN_{n,k}$.

The analysis of fixed points here will be similar to that in Subsection \ref{ca-pf-bouquet}, using vertex-orbits and edge-orbits of the action of elements of $\la S, R \ra$.  However, since edges now must join vertices of opposite sign (i.e., opposite parity) we have additional restrictions on the edge-orbits.

\begin{computing}{$\Aax(I)$}
To compute $\Aax(I)$ we consider fixed points of elements of $\Csa$, which have the form $S^h$.  If $h$ is odd, $S^h$ will move vertices to vertices of opposite sign, and hence will not preserve a triple $(N, \col, \ep)$.  Therefore, we can only have fixed points if $h$ is even. If $h$ is even we know signs are preserved, and there are $2$ choices of $\ep$, so we just need to determine when $(N, \col)$ is fixed.

Therefore, we may suppose that $h=2j$ is even.   The vertex-orbits of $\Ztn$ under the action of $\la S^h \ra$ are just the cosets $a + 2j\Ztn = a + 2g\Ztn$, where $g = (j, n)$.  Each coset has size $d = 2n/2g = n/g$.  There are $\eul(d)$ possible values of $j$ for each given $g$ and $d$ with $gd = n$.  Since $2g$ is even, all elements of each coset have the same sign, and there are $g$ even cosets and $g$ odd cosets.  We can therefore construct a fixed point of $S^h$ by matching even cosets to odd cosets in $g!$ ways, choosing an edge-orbit for each pair of cosets in $d^g$ ways, and choosing a coloring of the edge-orbits in $k^g$ ways.  Thus, the number of fixed $(N, \col)$ for $S^h$ is $g! d^g k^g$.

Adding over all possible even $h$, and remembering that $\ep$ can be chosen in $2$ ways, gives
$$\Aax(I) = \frac{1}{2n} \kern4pt 2 \kern-8pt 
		\sum_{(d,g): dg = n} \phi(d)\, g!\, d^g\, k^g
	= \frac{1}{n} \sum_{(d,g): dg = n} \phi(d)\, g!\, d^g\, k^g .$$
\end{computing}

\begin{computing}{$\Aax(R)$}
Suppose $\ga = S^h R \in \Csa R$ and $(N, \col, \ep) \in \fix(\ga)$. Since $R$ preserves the parity, hence the sign, of elements of $\Ztn$, again $h$ must be even to have a fixed point.  If $h$ is even we know $S^h R$ preserves signs, and there are $2$ choices of $\ep$, so we just need to determine when $(N, \col)$ is fixed.

Therefore, we may suppose that $h = 2j$ is even.  For $i \in \Ztn$, we have $S^hR(i) = h-i$, and the vertex-orbits of $S^h R$ are $\set{j}$ and $\set{j+n}$ of size $1$, and $n-1$ vertex-orbits of size $2$.  The two vertex-orbits of size $1$ must be matched by $N$ to each other.  This is only possible if one is even and the other odd, which means $n$ must be odd; if $n$ is even we have no fixed points.

Thus, we may assume $n$ is odd, so that there is one even and one odd vertex-orbit of size $1$, and $(n-1)/2$ even and $(n-1)/2$ odd vertex-orbits of size $2$. We can therefore construct a fixed point of $S^h R$ by matching the two vertex-orbits of size $1$ in a unique way, matching even vertex-orbits of size $2$ to odd vertex-orbits of size $2$ in $((n-1)/2)!$ ways, choosing edge-orbits for these pairs in $2^{(n-1)/2}$ ways, and then coloring the $1+(n-1)/2 = (n+1)/2$ edge-orbits in $k^{(n+1)/2}$ ways.

Adding over all $n$ possible even values of $h$, and remembering that $\ep$ can be chosen in $2$ ways, gives
$$\Aax(R) = \begin{cases}
 0 & \text{if $n$ is even,} \\
 \displaystyle \frac{1}{2n}\, 2n\, \left(\frac{n-1}{2}\right)!\,
2^{(n-1)/2}\,
		k^{(n+1)/2}
     = \left(\frac{n-1}{2}\right)!\, 2^{(n-1)/2}\, k^{(n+1)/2}
   & \text{if $n$ is odd.} \\
\end{cases}
$$
\end{computing}

\begin{computing}{$\Aax(F)$}
Suppose $\ga = S^h F \in \Csa F$ and $(N, \col, \ep) \in \fix(\ga)$.  Now $\ga$ flips signs, so we have the opposite situation to when $\ga = S^h$: we only obtain fixed points when $h$ is odd.  When $h$ is odd, we know $S^h F$ preserves signs, and there are $2$ choices of $\ep$, so we just need to determine when $(N, \col)$ is fixed.

Therefore, we may suppose that $h$ is odd.  The vertex-orbits of $\Ztn$ under the action of $\la S^h F\ra$ are just the cosets $a + h\Ztn = a+ g\Ztn$ where $g = (h,2n)$ is odd, so $g = (h, n) | n$.  Each coset has size $2d$ where $d = n/g$, and consists of $d$ even and $d$ odd elements.  For a given $g$ and $d$ there are $\eul(2d)$ possible values of $h$.

An edge-orbit of $N$ can match two cosets, or match a coset to itself.  Since there are $g$ cosets and $g$ is odd, at least one coset must be matched to itself. But if a coset is matched to itself, then each element $a'$ must be matched by $N$ to $a'+n$ (as in the analysis of $\Aax(I)$ in Subsection \ref{ca-pf-bouquet}).  Since $a'$ and $a'+n$ must have different signs, $n$ must be odd, and hence when $n$ is even there are no fixed points. If a coset $a+g\Ztn$ is matched to a different coset $b+g\Ztn$, then there are $d$ choices for the edge $\set{a,b+ug} \in N$, which determines the other edges in this edge-orbit.

Therefore, when $n$ is odd the $(N, \col)$ fixed by $\ga$ can be constructed as follows.  Choose $j$ with $0 \le j \le \lfloor g/2 \rfloor$, match $j$ pairs of cosets in $\ma(g, j)$ ways, choose the edge-orbits for these cosets in $d^j$ ways, then the other $g-2j$ edge-orbits which match cosets to themselves are uniquely determined, and the $j+(g-2j) = g-j$ edge-orbits can be colored in $k^{g-j}$ ways.  Thus, the number of fixed $(N, \col)$ for $S^h F$ is
$$\sum_{j=1}^{\lfloor g/2 \rfloor} \ma(g,j)\, d^j\, k^{g-j} .$$

Adding over all possible odd $h$, and remembering that $\ep$ can be chosen in $2$ ways, gives
$$
 \Aax(F) =
   \displaystyle \frac{1}{2n} \kern4pt
	2 \kern-8pt \sum_{\substack{(d,g)\,:\, dg = n \\ \text{$g$ odd}}}
	\eul(2d) \sum_{j=0}^{\lfloor g/2 \rfloor}
	\ma(g, j)\, d^j\, k^{g-j}
	 \quad\text{if $n$ is odd.}
$$
Simplifying gives the general expression
$$
 \Aax(F) 
  = \begin{cases}
	0 & \text{if $n$ is even,} \\
   \displaystyle \frac{1}{n}
	\sum_{\substack{(d,g)\,:\, dg = n \\ \text{$g$ odd}}}
	\eul(2d) \sum_{j=0}^{\lfloor g/2 \rfloor}
	\ma(g, j)\, d^j\, k^{g-j}
	& \text{if $n$ is odd.} \\
   \end{cases}
$$
\end{computing}

\begin{computing}{$\Aax(RF)$}
Suppose $\ga = S^h R F \in \Csa R F$ and $(N, \col, \ep) \in \fix(\ga)$. Since $\ga$ flips signs and $R$ preserves signs, we only obtain fixed points when $h$ is odd.  When $h$ is odd, we know $S^h R F$ preserves signs, and there are $2$ choices of $\ep$, so we just need to determine when $(N, \col)$ is fixed.

Therefore, we may suppose that $h$ is odd.  The vertex-orbits of $\Ztn$ under the action of $S^h R F$ are just the $n$ pairs $\set{i, h-i}$, each of which contains one even element and one odd element. An edge-orbit of $N$ can match two vertex-orbits, but also match a vertex-orbit to itself.

So $(N, \col)$ fixed by $\ga$ can be constructed as follows.  Choose $j$ with $0 \le j \le \lfloor n/2 \rfloor$, match $j$ pairs of vertex-orbits in $\ma(g, j)$ ways, and match the $n-2j$ remaining vertex-orbits to themselves.  We have only one choice for the edge-orbit for a pair of matched vertex-orbits or for the edge-orbit matching a vertex-orbit to itself.  We can color the $j + (n-2j) = n-j$ edge-orbits in $k^{n-j}$ ways. Thus, the number of $(N, \col)$ fixed by $S^h R F$ is
$$\sum_{j=1}^{\lfloor n/2 \rfloor} \ma(n,j)\, k^{n-j} .$$

Adding over the $n$ possible odd values of $h$, and remembering that $\ep$ can be chosen in $2$ ways, gives
$$\Aax(RF) = \frac{1}{2n} \kern4pt 2n
	\sum_{j=1}^{\lfloor n/2 \rfloor}
	\ma(n,j)\, k^{n-j}
 = \sum_{j=1}^{\lfloor n/2 \rfloor} \ma(n,j)\, k^{n-j} .$$
\end{computing}

\subsection{Counting directed embeddings of colored bouquets and related objects}\label{bo-pf-dirbouquet}

In this subsection we consider $\cA_{n,k}$, or equivalently $\cN_{n,k}$, under the action of a group $\Ga$ with $\Csa \le \Ga \le \Wga$.  In a similar way to previous sections, we can just take the average over all cosets of $\Csa$ in $\Ga$ of the average number of fixed points for a coset.

\begin{counting}{(A1) oriented colored directed embeddings of directed bouquets}
Given an oriented directed embedding $\Phi$ of a colored directed bouquet, we can transform it into a colored signed labeled bouquet by choosing a half-arc to label $0$, labeling the other half-arcs in ascendingi clockwise order, and then transforming the arc directions into positive signs on outward half-arcs and negative signs on inward half-arcs.  However, our choice of which half-arc to label $0$ was arbitrary.  So other colored signed labeled bouquets for $\Phi$ can be obtained by applying an arbitrary cyclic shift $S^h$.

Thus, an oriented colored directed embeddings of a directed bouquett may be regarded as an equivalence class of colored signed labeled bouquets under the action of $\Csa = \la S \ra$, and so the number of equivalence classes is just $\Aax(I)$.
\end{counting}

\begin{counting}{(A2) orientable colored directed embeddings of directed bouquets}
To count orientable directed embeddings, we must allow for the orientation of the surface being reversed.  As with (B2) earlier, we have to add the transformation $R$ to our group of symmetries, giving $\la S, R \ra = \Csa \cup \Csa R$.
The number of equivalence classes is therefore
$\rc2 (\Aax(I)+\Aax(R))$.
\end{counting}

\begin{counting}{(A3) arc reversal classes of oriented colored directed embeddings of directed bouquets}
To allow for arc reversal, we must add the transformation $F$ to our group of symmetries, giving $\la S, F \ra = \Csa \cup \Csa F$.
The number of equivalence classes is therefore
$\rc2 (\Aax(I)+\Aax(F))$.
\end{counting}

\begin{counting}{(A4) arc reversal classes of orientable colored directed embeddings of directed bouquets}
To allow for both surface orientation reversal and arc reversal, we must add both $R$ and $F$ to our group of symmetries, giving $\la S, R, F\ra = \Csa \cup \Csa R \cup \Csa F \cup \Csa RF$.
The number of equivalence classes is therefore
$\rc4 (\Aax(I) + \Aax(R) + \Aax(F) + \Aax(RF))$.
\end{counting}

\begin{counting}{(A5) simultaneous reflection and arc reversal classes of oriented colored directed embeddings of directed bouquets}
To allow for simultaneous reflection and arc reversal, we must add the transformation $RF$ to our group of symmetries, giving $\la S, RF \ra = \Csa \cup \Csa RF$.
The number of equivalence classes is therefore
$\rc2 (\Aax(I)+\Aax(RF))$.
\end{counting}

There are five groups $\Ga$ that satisfy $\Csa \le \Ga \le \Wga$.  These correspond to the subgroups of the quotient group $\Wga/\Csa$, which is a $4$-element dihedral group (or Klein group).  We have counted equivalence classes of $\cA_{n,k}$ (or $\cN_{n,k}$) under the action of $\Ga$ for all five groups in (A1)--(A5).

\begin{counting}{(A6) generic (orientable or nonorientable) colored directed embeddings of directed bouquets}
Similarly to (B3), generic $k$-colored directed embeddings of directed bouquets are in one-to-one correspondence with $2k$-colored orientable directed embeddings of directed bouquets.  So this is (A2) with $k$ replaced by $2k$, namely
$\rc2 ( \Aa(I)(n,2k)+\Aa(R)(n,2k) )$.
\end{counting}

\begin{counting}{(A7) nonorientable colored directed embeddings of directed bouquets}
Similarly to (B4), this is
$\text{(A6)} - \text{(A2)} =
	\rc2 ( \Aa(I)(n,2k)+\Aa(R)(n,2k) - \Aa(I)(n,k)-\Aa(R)(n,k) )$.
\end{counting}

\begin{counting}{(A8) arc-reversal classes of generic (orientable or
nonorientable) colored directed embeddings of directed bouquets}
Similarly to (B3) and (A6), generic $k$-colored arc-reversible directed embeddings of directed bouquets are in one-to-one correspondence with the corresponding $2k$-colored orientable objects.  So this is (A4) with $k$ replaced by $2k$, namely
$\rc4 (\Aa(I)(n,2k) + \Aa(R)(n,2k) + \Aa(F)(n,2k) + \Aa(RF)(n,2k))$.
\end{counting}

\begin{counting}{(A9) arc-reversal classes of nonorientable colored directed embeddings of directed bouquets}
Similarly to (B4) and (A7), this is
$\text{(A8)} - \text{(A4)} =
 \rc4 (\Aa(I)(n,2k) + \Aa(R)(n,2k) + \Aa(F)(n,2k) + \Aa(RF)(n,2k)
  -\Aa(I)(n,k) - \Aa(R)(n,k) - \Aa(F)(n,k) - \Aa(RF)(n,k))$.
\end{counting}

\begin{observation}
From (A1), (A2), (A3), and (A5) we see that $\Aax(I)$, $\Aax(R)$, $\Aax(F)$, and $\Aax(RF)$ are all integers.
\end{observation}

\printbibliography

\end{document}